\documentclass[final]{siamltex}
\usepackage{xcolor}
\usepackage{amsmath,amssymb,amsfonts}
\usepackage{mathrsfs}
\usepackage{graphicx}
\usepackage{bm}
\usepackage{ucs}
\usepackage[utf8x]{inputenc}
\usepackage{wrapfig}
\usepackage{color}
\usepackage{float}
\usepackage{epstopdf}
\usepackage{algorithm,algorithmic}
\usepackage[breaklinks,colorlinks=true,linkcolor=blue,citecolor=red,
  backref=page]{hyperref}

\DeclareMathAlphabet{\itbf}{OML}{cmm}{b}{it}

\newcommand{\RR}{\mathbb{R}}
\newcommand{\CC}{\mathbb{C}}
\newcommand{\PP}{\mathbb{P}}
\newcommand{\EE}{\mathbb{E}}

\newcommand{\NN}{\mathbb{N}}

\newcommand{\eps}{\varepsilon}

\newtheorem{remark}{Remark}[section]

\newcommand\balpha{{\boldsymbol{\alpha}}}

\begin{document}
\title{Wave propagation in randomly perturbed weakly coupled waveguides} \author{Liliana Borcea\footnotemark[1]
  \and Josselin
  Garnier\footnotemark[2]}

\maketitle


\renewcommand{\thefootnote}{\fnsymbol{footnote}}
\footnotetext[1]{Department of Mathematics, University of Michigan,
  Ann Arbor, MI 48109. {\tt borcea@umich.edu}}
\footnotetext[2]{Centre de Math\'ematiques Appliqu\'ees, Ecole Polytechnique, 91128 Palaiseau Cedex, France.  {\tt
    josselin.garnier@polytechnique.edu}}
\markboth{L. Borcea  and J. Garnier}{Coupled waveguides}

\maketitle

\begin{abstract}
We present an analysis of wave propagation in a two step-index, parallel waveguide system.  The goal  is to quantify the effect of scattering at randomly perturbed interfaces   
between the guiding layers of high index of refraction and the host medium.  
The analysis  is based on the expansion of the solution of the wave equation in a complete set of guided, radiation and evanescent modes with amplitudes that are random fields, due to scattering. We obtain a detailed characterization of these amplitudes and thus 
quantify the transfer of power between the two waveguides in terms of their separation distance. The results show that, 
no matter how small the  fluctuations of the interfaces are, they have  significant effect at sufficiently large distance of propagation,  which manifests in two ways: The first effect is well known and consists of power leakage from the guided modes to the radiation ones. The second effect  consists of blurring of the periodic transfer of power between the waveguides and the eventual equipartition of power. 
Its quantification is the main practical result of the paper.
\end{abstract}


\section{Introduction}
\label{sec:intro}
Guided waves have applications in electromagnetics
\cite{collin1990field}, optics and communications \cite{marcuse,syms}, 
imaging underwater \cite{kohler77,BorGar,gomez09,gomez11}, imaging of and in tunnels \cite{bedford2014modeling} and so on. 
The classical theory of guided waves is for ideal waveguides with perfectly reflecting straight walls and filled with 
homogeneous media, where the wave equation can be solved using separation of variables. The wave field is represented as a superposition of finitely many guided modes, which are waves that propagate along the axis of the waveguide, and infinitely many evanescent modes which decay away from the source. These modes do not interact with each other and thus have constant amplitude determined by the wave source \cite{collin1990field,syms}. 

Motivated by applications in imaging and communications, the classical theory has been extended to waveguides filled with random media \cite{marcuse,garnier07a,garnier07b,AlBor}, with randomly perturbed boundaries \cite{AlBorGar,AlBor} and with slowly changing cross-section \cite{BorGarTP,BorGarWTP,marcuse}. Weakly guiding
waveguides with confining graded-index profile affected by small random perturbations have also been analyzed in  \cite{garnier00,perrey07}.
The resulting mode coupling theory  quantifies the interaction between the  modes induced by scattering, and the consequent randomization of the wave field.  

We consider waveguides with penetrable boundaries, where the guiding effect is due to a medium of 
high index of refraction embedded in a homogeneous background. Such waveguides are analyzed in \cite{kohler77,gomez11,gomezRB} in the context of underwater acoustics \cite{wilcox} and are of great importance in optics and communications \cite{marcuse,Rowe99,syms}. Motivated 
by the latter applications, we consider a waveguide system made of two parallel step-index waveguides, which is known as a directional coupler 
in integrated optics \cite{Microring}, \cite[Chapter 10]{syms}. The classical analysis of this system, described in  \cite{huang94}, \cite[Chapter 10]{syms} and 
\cite[Chapter 10]{marcuse}, is based on the observation that the transverse profiles of the guided waves are essentially supported in the step-index waveguides and decay outside.  When the step-index waveguides are nearby, the decaying tails penetrate the 
neighboring waveguide and transfer of power can occur. This is the sole coupling mechanism in ideal directional couplers and for 
synchronous waveguides (with identical guided mode phase velocity) there is a complete, periodic transfer of power.  That is to say, 
if the source emits power in one step-index waveguide, this is transferred to the other waveguide and back in a periodic manner, at regular distance intervals.  These intervals depend on the separation between the step-index waveguides. The larger this distance is, the weaker the coupling and the farther the waves must travel for the transfer of power to occur.

We introduce a mathematical analysis of a randomly perturbed directional coupler, where the interfaces that 
separate the medium with high index of refraction  from the background have small amplitude random fluctuations on a scale 
similar to the wavelength.  The classic approach in \cite[Chapter 10]{syms}, which is based solely on the guided modes, is inadequate in this case, because scattering at the random interfaces induces mode coupling. We take into account all the modes, 
the guided, radiation and evanescent ones,  and quantify how their interaction affects  the performance of the directional coupler. 
The analysis is focused on the case of well separated waveguides, where the deterministic coupling is weak. It applies to an arbitrary 
number of guided modes, but we describe in depth the results for the case of single guided mode step-index waveguides. 
We show that mode coupling induced by the random fluctuations is present no matter how far apart the waveguides are and it has two effects:
The first effect is well known \cite[Chapter 10]{marcuse} and consists of power leakage
from the guided modes to the radiation modes. Our analysis captures it  and shows that the leaked power is self-averaging
i.e., it is independent of the realization of the random processes that model the fluctuations of the interfaces. The other effect consists 
of the blurring of the periodic transfer of power between the waveguides and the eventual equipartition of power between the guided waves. Its quantification in terms of the waveguide separation and amplitude of the random fluctuations is the main practical result of the paper. 

The paper is organized as follows: We begin in section \ref{sect:form} with the mathematical formulation of the problem 
and then state the results in section \ref{sect:results}. Their derivation is in sections \ref{sect:homog}--\ref{sect:randomTwo}. 
We end with a summary in section \ref{sect:sum}.

\section{Formulation of the problem}
\label{sect:form}

\begin{figure}[t]
\vspace{-3.5in}
\begin{picture}(0,0)%
\hspace{1.5in}\includegraphics[width = 0.5\textwidth]{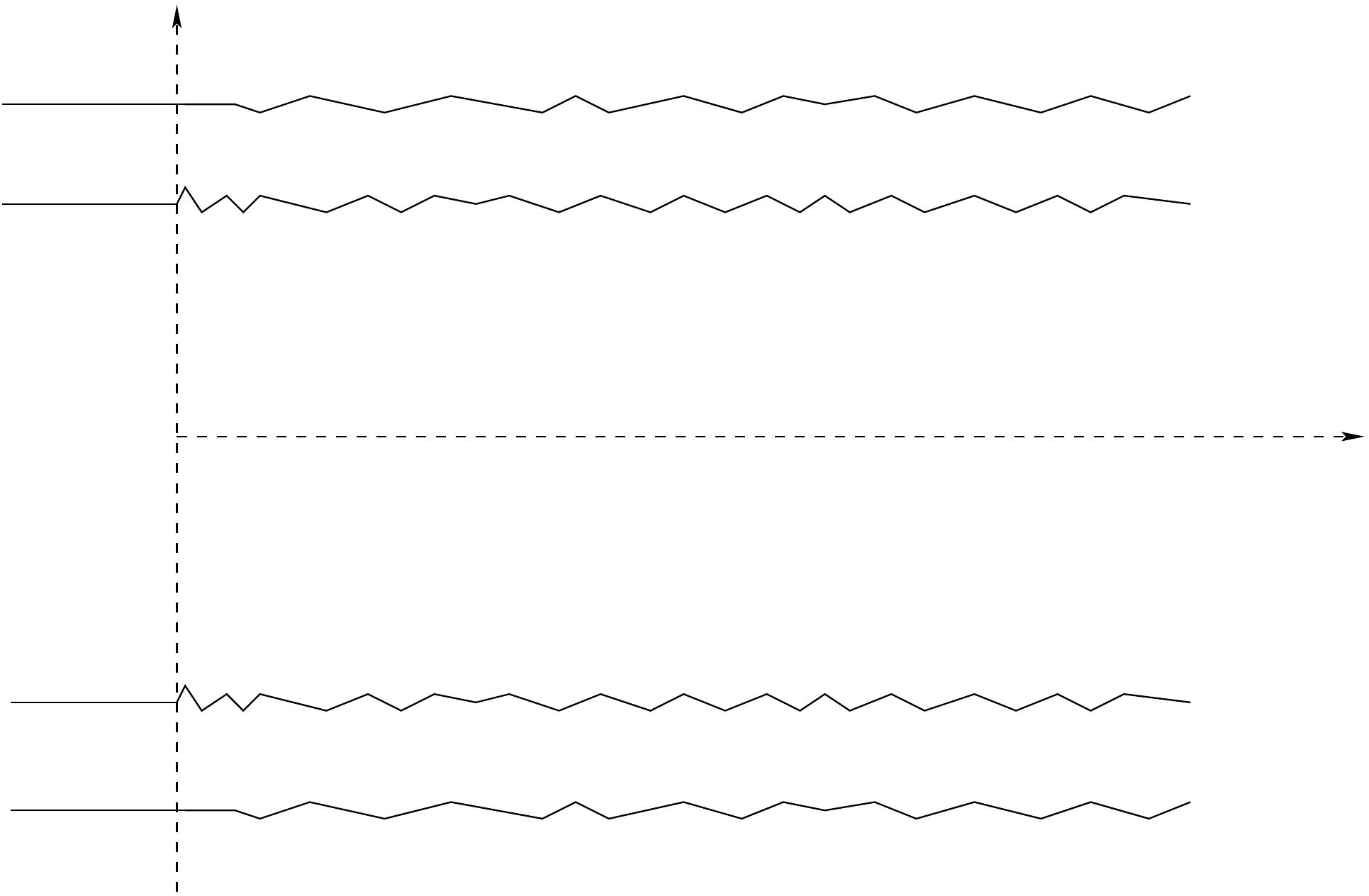}%
\end{picture}%
\setlength{\unitlength}{2960sp}%
\begingroup\makeatletter\ifx\SetFigFont\undefined%
\gdef\SetFigFont#1#2#3#4#5{%
  \reset@font\fontsize{#1}{#2pt}%
  \fontfamily{#3}\fontseries{#4}\fontshape{#5}%
  \selectfont}%
\fi\endgroup%
\begin{picture}(10799,8069)(1599,-7583)
\put(3300,-5200){\makebox(0,0)[lb]{\smash{{\SetFigFont{7}{8.4}{\familydefault}{\mddefault}{\updefault}{\color[rgb]{0,0,0}{\normalsize $\frac{d}{2} + D$}}%
}}}}
\put(3700,-5600){\makebox(0,0)[lb]{\smash{{\SetFigFont{7}{8.4}{\familydefault}{\mddefault}{\updefault}{\color[rgb]{0,0,0}{\normalsize $\frac{d}{2} $}}%
}}}}
\put(3110,-7400){\makebox(0,0)[lb]{\smash{{\SetFigFont{7}{8.4}{\familydefault}{\mddefault}{\updefault}{\color[rgb]{0,0,0}{\normalsize $-\frac{d}{2} - D$}}%
}}}}
\put(3600,-7100){\makebox(0,0)[lb]{\smash{{\SetFigFont{7}{8.4}{\familydefault}{\mddefault}{\updefault}{\color[rgb]{0,0,0}{\normalsize $-\frac{d}{2} $}}%
}}}}
\put(4300,-5000){\makebox(0,0)[lb]{\smash{{\SetFigFont{7}{8.4}{\familydefault}{\mddefault}{\updefault}{\color[rgb]{0,0,0}{\normalsize $x$}}%
}}}}
\put(8000,-6400){\makebox(0,0)[lb]{\smash{{\SetFigFont{7}{8.4}{\familydefault}{\mddefault}{\updefault}{\color[rgb]{0,0,0}{\normalsize $z$}}%
}}}}
\put(5500,-5400){\makebox(0,0)[lb]{\smash{{\SetFigFont{7}{8.4}{\familydefault}{\mddefault}{\updefault}{\color[rgb]{0,0,0}{\normalsize ${\rm n}^{(\eps)}(z,x) = n$}}%
}}}}
\put(5500,-7200){\makebox(0,0)[lb]{\smash{{\SetFigFont{7}{8.4}{\familydefault}{\mddefault}{\updefault}{\color[rgb]{0,0,0}{\normalsize ${\rm n}^{(\eps)}(z,x) = n$}}%
}}}}
\put(5500,-7600){\makebox(0,0)[lb]{\smash{{\SetFigFont{7}{8.4}{\familydefault}{\mddefault}{\updefault}{\color[rgb]{0,0,0}{\normalsize ${\rm n}^{(\eps)}(z,x) = 1$}}%
}}}}
\put(5500,-6000){\makebox(0,0)[lb]{\smash{{\SetFigFont{7}{8.4}{\familydefault}{\mddefault}{\updefault}{\color[rgb]{0,0,0}{\normalsize ${\rm n}^{(\eps)}(z,x) = 1$}}%
}}}}
\put(5500,-5000){\makebox(0,0)[lb]{\smash{{\SetFigFont{7}{8.4}{\familydefault}{\mddefault}{\updefault}{\color[rgb]{0,0,0}{\normalsize ${\rm n}^{(\eps)}(z,x) = 1$}}%
}}}}
\end{picture}%
\vspace{0.1in}
\caption{Illustration of two waveguides with fluctuating interfaces,
  filled with a medium with index of refraction $n > 1$. The waves
  propagate along the range axis $z$. The waveguides have width $D$ and are
  separated by the distance $d$.  }
\label{fig:setup}
\end{figure}
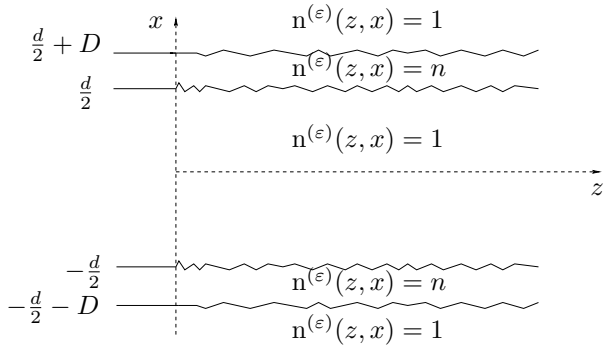

We study the propagation of a time harmonic wave in a medium with
index of refraction ${\rm n}^{(\eps)}(z,x)$ defined below, that models two step-index 
waveguides of width $D$, separated in the transverse direction $x$ by
the distance $d$, as illustrated in Figure \ref{fig:setup}. The wave
field is denoted by $p(z,x)$ and it solves the two-dimensional Helmholtz equation
\begin{align}
\label{eq:fouriertransform}
\big(\partial_x^2+\partial_z^2\big){p}(z,x) + \big[k{\rm n}^{(\eps)}(z,x)\big]^2 {p}(z,x) =
    {f}(x) \delta'(z) ,
\end{align}
with radiation condition at infinity, where $k$ is the wavenumber
 and $f(x)$ models a source supported
at the origin of the range coordinate $z$.

In the case of ideal (unperturbed) waveguides, the index of refraction
is range independent and equal to
\begin{align}
{\rm n}^{(0)}(x) = \left\{
\begin{array}{ll}
\it{n} & \mbox{ if } x \in (-d/2-D,-d/2)\cup (d/2,d/2+D),\\
1 &\mbox{ otherwise,}
\end{array}
\right.
\label{eq:no}
\end{align}
with $n>1$. We are interested in perturbed waveguides, where the index
of refraction
\begin{equation}
\label{eq:modelpert1}
{\rm n}^{(\eps)}(z,x) = \left\{
\begin{array}{ll}
n & \mbox{ if } x \in ({\cal D}_{1}^{(\eps)}(z), {\cal D}_{2}^{(\eps)}(z) )\cup ({\cal
  D}_{3}^{\eps}(z), {\cal D}_{4}^{(\eps)}(z) ) , \\ 1 &\mbox{ otherwise,}
\end{array}
\right.
\end{equation}
jumps across four randomly fluctuating interfaces 
\begin{align}
{\cal D}_{1}^{(\eps)}(z) &= -d/2-D+ \eps D\nu_1(z) {\bf 1}_{(0,L^{(\eps)})}(z), \nonumber \\
{\cal D}_{2}^{(\eps)}(z) &= -d/2  +\eps D\nu_2(z) {\bf 1}_{(0,L^{(\eps)})}(z),\nonumber \\
{\cal D}_{3}^{(\eps)}(z) &= d/2+\eps D\nu_3(z) {\bf 1}_{(0,L^{(\eps)})}(z) , \nonumber \\
{\cal D}_{4}^{(\eps)}(z) &= d/2 +D +\eps D\nu_4(z) {\bf 1}_{(0,L^{(\eps)})}(z) . \label{eq:Interfaces}
\end{align}
The fluctuations are modeled by the zero-mean, bounded,
independent and identically distributed stationary random processes $\{\nu_q(z)\}_{1 \le q \le 4}$ 
with smooth covariance function
\begin{equation}
{\cal R}(z) = \EE[\nu_q(0)\nu_q(z)], \qquad q = 1, \ldots, 4. \label{eq:covar}
\end{equation}
These satisfy strong mixing conditions as defined for example in
\cite[section 2]{papa74}. The typical amplitude of the fluctuations  is much smaller than $D$ and it is modeled in
\eqref{eq:Interfaces} by the
small and positive dimensionless parameter $\eps$.

We study the wavefield at $z > 0$, satisfying  
\begin{equation}
{p}(z,x)
\in {\cal C}^0\big( (0,+\infty),H_0^1(\RR)\cap H^2(\RR)\big)
\cap  {\cal C}^2\big( (0,+\infty),L^2(\RR)\big) , \qquad z > 0,
\label{eq:functclass}
\end{equation}
and to set radiation conditions, we suppose that the random
fluctuations are supported in the range interval $(0,L^{(\eps)})$.  We
will see that net scattering effect of these fluctuations  becomes of
order one at range scales of order $\eps^{-2}$, so we let $L^{(\eps)} =
L/\eps^2$. We will also see that for the assumed 
smooth covariance ${\cal R}(z)$, the guided
waves propagate mostly in the forward direction and do not interact
with any fluctuations at $z < 0$, which is why we neglect them.

The goal of the paper is to quantify how scattering at the random
interfaces \eqref{eq:Interfaces} affects the coupling of the 
two step-index waveguides centered at $x = \pm (d+D)/2$. We consider in particular the case of a
sufficiently large separation distance $d$ between the waveguides, where the
deterministic coupling is very weak but the coupling induced by the
random fluctuations is still present.

\section{Statement of results}
\label{sect:results}
We state here the main results of the paper, derived in sections
\ref{sect:homog}--\ref{sect:randomTwo} by decomposing the wavefield
into guided, radiation and evanescent modes of the waveguide system made of two parallel step-index waveguides 
illustrated in Figure \ref{fig:setup}. 
While our analysis applies to an arbitrary number of guided modes, in this section we consider
only the case where an isolated step-index waveguide 
has only one guided mode.  This case captures all the essential
aspects of the problem and arises when 
\begin{equation}
\label{eq:assumeD}
k D \sqrt{n^2-1} < \pi .
\end{equation}

The two ideal step-index waveguides centered at $\pm(d/2+D/2)$ do not interact in the waveguide system when $d \to \infty$. Thus,  we
can write the wavefield for large $d$ in terms of the unique guided mode $\phi(x) e^{i \beta z}$ 
of the step-index waveguide centered at $x=d/2+D/2$, modeled by the index of refraction 
\begin{align}
{{\rm n}}_s(x) = \left\{
\begin{array}{ll}
\it{n} & \mbox{ if } x \in (d/2,d/2+D),\\ 1 &\mbox{ otherwise,}
\end{array}
\right.
\label{eq:notilde}
\end{align}
and in terms of the  unique guided mode $\phi(-x) e^{i \beta z}$ of the step-index waveguide centered at $x=-d/2-D/2$, modeled by the index of refraction ${{\rm n}}_s(-x)$. 
Here $\phi(x)$ is the eigenfunction 
of the Helmholtz operator $\partial_x^2 + {\rm n}_s(x)$ for the  eigenvalue $\beta^2$. This $\beta$ is defined in  Lemma \ref{lem.singlebeta}
as the unique solution  in the interval $(k,nk)$ of 
\begin{equation}
\frac{\sqrt{n^2 k^2 - \beta^2}}{\sqrt{\beta^2 - k^2}} \tan \Bigg(
\frac{D}{2} \sqrt{k^2 n^2 - \beta^2}\Bigg) = 1, \qquad \beta \in (k,nk).
\label{eq:defbeta}
\end{equation}
The expression of the eigenfunction $\phi(x)$ is 
\cite[Section 2]{magnanini}, \cite[Chapter 1]{marcuse} 
\begin{align}
\phi(x) = 
\left\{
\begin{array}{ll}
\Big(\frac{2}{\eta} + D \Big)^{-1/2} \cos (\xi \frac{D}{2}) \exp [\eta
  (x-\frac{d}{2})] , \qquad  & x \le  \frac{d}{2} ,
\\ \Big(\frac{2}{\eta} + D \Big)^{-1/2} \cos
\big[\xi(x-\frac{d}{2}-\frac{D}{2})\big],  & \frac{d}{2}\leq
x\leq \frac{d}{2}+D ,\\ \Big(\frac{2}{\eta} + D \Big)^{-1/2} \cos (\xi
\frac{D}{2}) \exp\big[- \eta (x-\frac{d}{2}-D)\big], & x\geq
\frac{d}{2}+D.
\end{array}
\right.
\label{eq:phi}
\end{align}
It has  a peak centered at $x=d/2+D/2$ of width  $1/\xi$, and an exponentially decaying tail, at the rate $\eta$, where 
\begin{equation}
\xi = \sqrt{n^2 k^2 - \beta^2}, \quad \eta = \sqrt{\beta^2 - k^2}.
\label{eq:defxieta}
\end{equation}

\subsection{Coupling of ideal waveguides}
\label{sect:resHomog}
\begin{figure}
\begin{center}
\includegraphics[width=5.5cm]{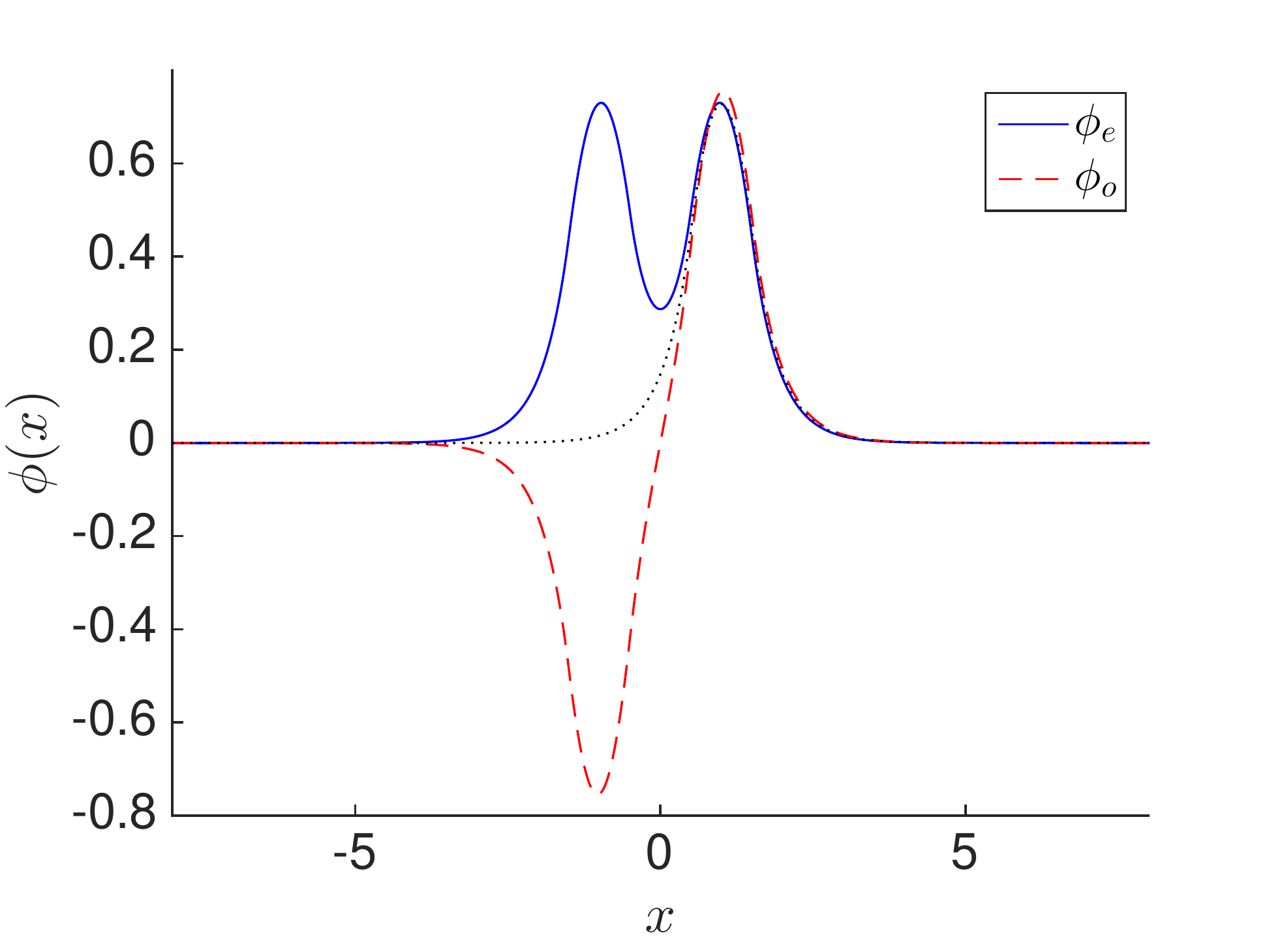} 
\includegraphics[width=5.5cm]{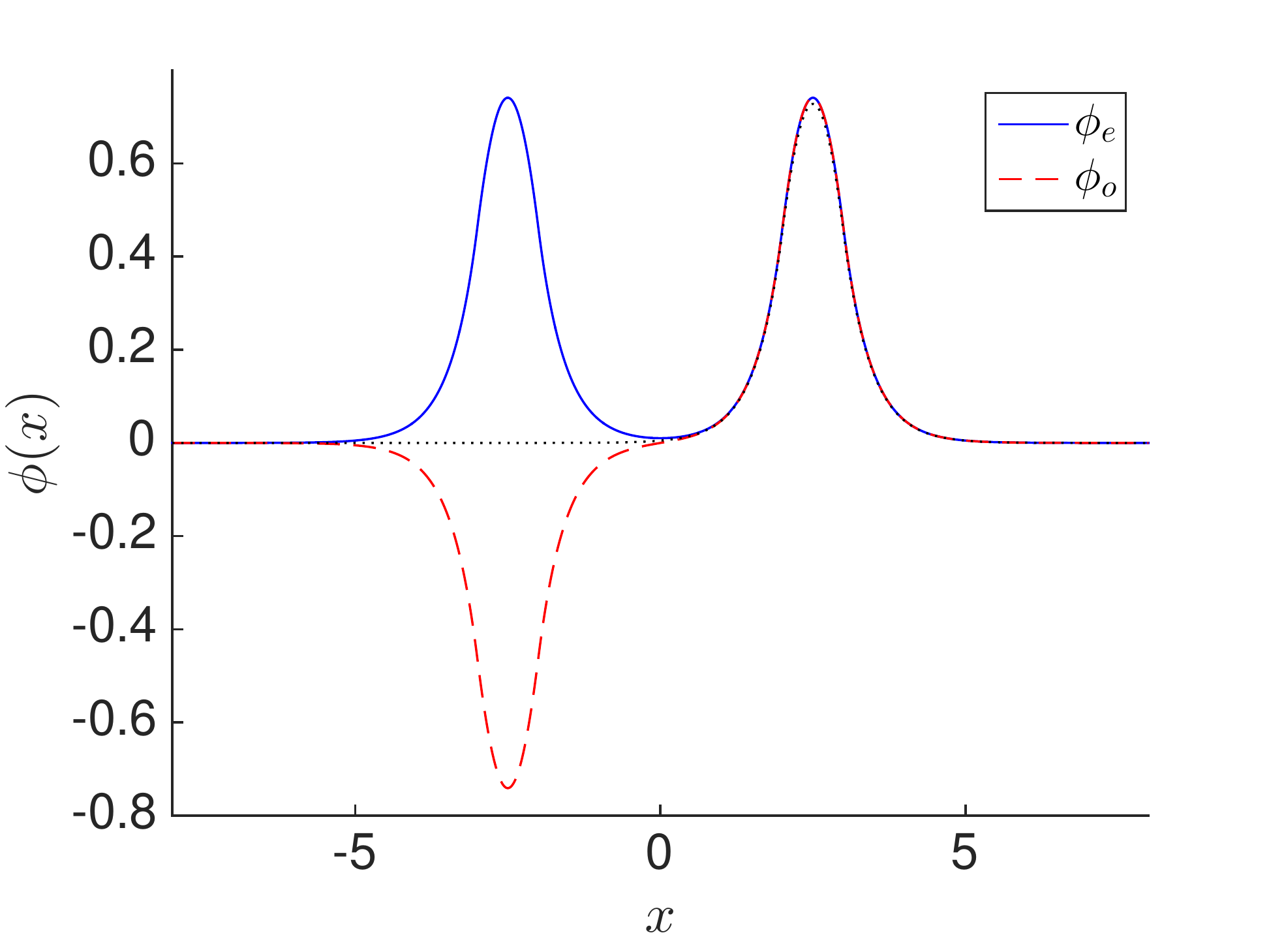} 
\end{center}
\caption{The eigenfunctions $\phi_e(x)$ (solid blue line) and $\phi_o(x)$ (dashed red  line) calculated 
at wavenumber $k = 2 \pi$, for the waveguide system with index of refraction $n = 1.1$. For reference we also plot 
the eigenfunction $\phi(x)$ of the step-index waveguide centered at $x=d/2+D/2$ 
with the black dotted line. The abscissa is in units of the waveguide width $D$, 
which is equal to the wavelength. We consider a separation $d = D$ (left plot) and $d= 4D$ (right plot)
between the  waveguides.}
\label{fig:eigenfunctions}
\end{figure}

We show in section \ref{sect:homog} that under the assumption \eqref{eq:assumeD} and for sufficiently large 
separation distance $d$, the solution of 
\begin{align}
\label{eq:homogHelm}
\big(\partial_x^2+\partial_z^2)p^{(0)}(z,x)+\big[k {\rm n}^{(0)}(x)\big]^2 {p}^{(0)}(z,x) =
    {f}(x) \delta'(z) ,
\end{align}
with radiation condition at infinity, takes the form\footnote{We use consistently the index ${(0)}$ for the mode amplitudes and wave field in the ideal (unperturbed) two step-index wave guide system.}
\begin{equation}
{p}^{(0)}(z,x) = \sum_{t \in \{e,o\}} \frac{{a}_{t}^{(0)}}{\sqrt{\beta_t}}
e^{i \beta_t z} \phi_t(x) +O(z^{-2}), \qquad z > 0.
\label{eq:phom}
\end{equation}
The terms in the sum model the guided modes of the waveguide system, 
whereas the $ O(z^{-2})$ remainder  accounts for the 
radiation and evanescent modes \cite[Section 3]{magnanini}.  
The guided modes  $\{\phi_t(x)e^{i \beta_t z}\}_{t = e,o}$ are waves that propagate along the range coordinate 
$z$ and have transverse profiles essentially supported in the two step-index waveguides. They are defined by  the even and odd eigenfunctions 
 $\phi_e(x)$ and $\phi_o(x)$  of the operator
\begin{align}
{\cal{H}}(x) = \partial_x^2 + {\rm n}^{(0)}(x), 
\label{eq:HelmIdeal}
\end{align}
with ${\rm n}^{(0)}(x)$ given in \eqref{eq:no}, for the eigenvalues
$\beta_e^2$ and $\beta_o^2$ satisfying
\begin{equation}
\beta_t \in (k, nk), \qquad t \in \{e,o\}.
\end{equation}
These eigenfunctions are given explicitly in section \ref{sect:discrSpect}. They
have peaks close to the center axes $x = \pm (d/2+D/2)$ of the step-index waveguides 
and decay outside their transverse support, as illustrated 
in Figure \ref{fig:eigenfunctions}. Each guided mode is multiplied in \eqref{eq:phom}
by the constants 
\begin{equation}
a_{t}^{(0)} = \frac{\sqrt{\beta_t}}{2} \int_{\RR} d x \,
\overline{\phi_t(x)} f(x), \quad t \in \{e,o\},
\label{eq:Amplo}
\end{equation}
called the mode amplitudes, that are determined by the source. Throughout the paper the bar denotes complex conjugate.  

We assume that  the waveguides separation distance $d$  is sufficiently large, so that 
\begin{equation}
\exp(-\eta d) \ll 1,
\label{eq:dlarge}
\end{equation}
with $\eta$ defined in \eqref{eq:defxieta}, and obtain that  the eigenfunctions can be approximated by 
\begin{align}
\phi_e(x) &= \phi(|x|) + O\big(e^{-\eta d}\big), \label{eq:phieoapr}
\quad \phi_o(x) = {\rm sgn}(x)\phi(|x|) + O\big(e^{-\eta d}\big),
\qquad x \in \RR,
\end{align}
where ``sgn" is the sign function. The accuracy of this approximation is illustrated   
in the right plot of Figure \ref{fig:eigenfunctions}. Consequently, 
the transverse profile of the even guided mode presents two
positive peaks centered at $x=\pm ( d/2+D/2)$, with exponentially
decaying tails, while the transverse profile of the odd guided mode
presents one positive peak centered at $x=d/2+D/2$, one negative
peak centered at $x=-d/2-D/2$, and exponentially decaying tails. The form of these peaks is proportional
to the unique eigenfunction \eqref{eq:phi} of the single-mode step-index
waveguide. 

We also obtain that  the wave numbers are
\begin{align}
\beta_{e} = \beta+ \beta' e^{-\eta d } +o\big( e^{-\eta d
}\big), \label{eq:betaeoapr} \quad \beta_{o} = \beta - \beta' e^{-\eta
  d }+o\big( e^{-\eta d }\big) , \end{align} with
\begin{align}
 \beta' &= \frac{\eta}{\beta \big(1+\frac{\eta^2}{\xi^2}\big)
   \big(\frac{1}{\eta}+\frac{D}{2}\big)}, \label{eq:Betaprime}
\end{align}
so  the wavefield \eqref{eq:phom} in the ideal waveguide system takes the form
\begin{equation}
\hspace{-0.07in}p^{(0)}(z,x) = \frac{\phi(|x|)}{\sqrt{\beta}} e^{i \beta z} \left[
  1_{(0,\infty)}(x) u_+^{(0)}(z) + 1_{(-\infty,0)}(x) u_{-}^{(0)}(z) \right] +
O\big(e^{-\eta d}\big) + O(z^{-2}).
\label{eq:Po1}
\end{equation}
Here we introduced the range-dependent amplitudes of the waves propagating in the two step-index waveguides
\begin{align}
u_+^{(0)}(z) &= (a_{e}^{(0)} + a_{o}^{(0)})\cos\big(\beta' z e^{-\eta d} \big) + i
(a_{e} ^{(0)}-a_{o}^{(0)})\sin\big(\beta' z e^{-\eta d} \big), \label{eq:up0}
\\ u_-^{(0)}(z) &= (a_{e}^{(0)} - a_{o}^{(0)})\cos\big(\beta' z e^{-\eta d} \big) +
i (a_{e}^{(0)} +a_{o}^{(0)})\sin\big(\beta' z e^{-\eta d} \big)\label{eq:um0},
\end{align}
with indexes ``$\pm$" corresponding
to the waveguide centered at $x = \pm (d+D)/2$.

Equation \eqref{eq:Po1} shows that 
the wavefield consists of two components, the first one is centered at $x=d/2+D/2$ with the form
$\phi(|x|) e^{i\beta z}{\bf 1}_{(0,\infty)}(x) $ and the second one is centered at $x=-d/2-D/2$ with the form
$\phi(|x|) e^{i\beta z}{\bf 1}_{(-\infty,0)}(x)$. This is similar to the case  of two independent, single-mode step-index waveguides,
except that in \eqref{eq:Po1} the amplitudes $u_\pm(z)$ vary in $z$, due to coupling.
We obtain from \eqref{eq:dlarge} and (\ref{eq:up0}--\ref{eq:um0}) that for $z$ of the order of the wavelength, 
\begin{align}
u_\pm^{(0)}(z) &\approx u_\pm^{(0)}(0) =  a_{e}^{(0)} \pm a_o^{(0)}.
\label{eq:uinit}
\end{align}
However,  at large  $z$, satisfying 
\begin{equation}
z  = \exp(\eta d) Z, \quad Z > 0, 
\label{eq:largez}
\end{equation}
$u_{\pm}^{(0)}(z)$ oscillate periodically in $z$. For example, if the
source gives the amplitudes 
\begin{equation}
a_{e}^{(0)} = a_{o}^{(0)} = \frac{a^{(0)}}{2}, \label{eq:assumea0}
\end{equation}
according to \eqref{eq:Amplo},
so that 
$u_+^{(0)}(0) = a^{(0)}$ and $u_-^{(0)}(0)  = 0$, 
at range  \eqref{eq:largez} we have
\begin{align}
u_+^{(0)}\big( e^{\eta d} Z \big) = a^{(0)}\cos(\beta' Z), \quad \quad
u_-^{(0)}\big( e^{\eta d} Z \big) =  i a^{(0)} \sin(\beta' Z).
\label{eq:uoscZhom}
\end{align}

In conclusion, the total wave power in the ideal waveguide system at large range $z$ is essentially supported in the two step-index waveguides.  The wave power  in the waveguide centered at $x = \pm(d/2+D/2)$ is proportional to $|u_\pm^{(0)}(z)|^2$, and 
equation  (\ref{eq:uoscZhom}) shows that  it oscillates slowly and periodically.  At scaled distance $Z =m \pi / \beta'$, $m\in \NN$, 
the wave power is concentrated in the
waveguide centered at $x=d/2+D/2$, whereas at $Z = (1/2+m) \pi /
\beta '$, $m\in \NN$, the wave power is concentrated
in the waveguide centered at $x=-d/2-D/2$.  These periodic oscillations
have been reported in the literature \cite[Chapter 10]{syms}.  The
standard method to analyze them is not to start from the analysis of
the modes of the waveguide system, as we do in section
\ref{sect:homog}, but to simplify by assuming that the modes  can be represented as a weighted sum of the guided
modes of the two  waveguides.  This simplified approach does
not allow to take into account the role of evanescent and radiation
modes, which are critical to the study of random waveguides in sections
\ref{sect:random}--\ref{sect:randomTwo}, with results described next.

\subsection{Coupling of random  waveguides}
\label{sect:resRand}
The analysis of the solution $p(z,x)$ of the Helmholtz equation
\eqref{eq:fouriertransform} with index of refraction
\eqref{eq:modelpert1} is carried out in sections \ref{sect:random}--\ref{sect:randomTwo}. It
shows that under the assumptions \eqref{eq:assumeD}, \eqref{eq:dlarge}
and at large range $z/\eps^2$, where scattering at the random
interfaces \eqref{eq:Interfaces} becomes significant, there are three
distinguished regimes that determine the coupling between the 
waveguides:
\vspace{0.05in}
\begin{enumerate}
\item The ``moderate coupling" regime, where the separation distance $d$  is moderately large, satisfying
\begin{equation}
1 \gg \exp(-\eta d) \gg \eps^2.
\label{eq:modReg}
\end{equation}
\item The  ``weak coupling" regime, where $d$ is large enough so that 
\begin{equation}
\exp(-\eta d) = O(\eps^2).
\label{eq:larged}
\end{equation}
\item The  ``very weak coupling" regime, where $d$ is so large that 
\begin{equation}
\exp(-\eta d) \ll \eps^2 \ll 1.
\label{eq:huged}
\end{equation}
\end{enumerate}
We now describe the results in each of these three regimes.

\subsubsection{Moderate coupling} 
\label{sect:resModer}
At large range $z/\eps^2$ and in the regime defined by \eqref{eq:modReg}, 
the solution of \eqref{eq:fouriertransform} with radiation condition at infinity and with index of refraction  \eqref{eq:modelpert1}
is 
\begin{equation}
p\Big(\frac{z}{\eps^2},x\Big) = \frac{\phi(|x|)}{\sqrt{\beta}} e^{i
  \beta \frac{z}{\eps^2}} \big[ 1_{(0,\infty)}(x) u_+(z) +
  1_{(-\infty,0)}(x) u_-(z) \big] + o(1),
\label{eq:pmoder}
\end{equation}
where  $u_\pm(z)$ are random processes and $o(1)$ denotes a residual
that tends to zero as $\eps \to 0$.  This residual accounts for the
radiation and evanescent components of $p$.

Similar to \eqref{eq:Po1}, we have a propagating wave
$\phi(|x|)e^{i \beta z/\eps^2}{\bf 1}_{(0,\infty)}(x)$ centered at $x = d/2+D/2$ and 
 another wave
$\phi(|x|)e^{i \beta z/\eps^2}{\bf 1}_{(-\infty,0)}(x)$ centered at $x=-d/2-D/2$.
The coupling between the two waveguides is described by 
the random variations of the complex amplitudes $u_\pm(z)$, the analogues of 
(\ref{eq:up0}--\ref{eq:um0}).  To write their expressions,
we introduce the notation
\begin{equation}
\Delta \beta_t = \beta_t - \beta, \qquad t \in \{e,o\},
\label{eq:DeltaBeta}
\end{equation}
which takes into account that the residual in \eqref{eq:betaeoapr} is
not negligible under the assumption \eqref{eq:modReg} at range $O(\eps^{-2})$. We have
\begin{align}
u_\pm(z) &= a_e(z) e^{i \Delta \beta_e \frac{z}{\eps^2}} \pm
a_o(z) e^{i \Delta \beta_o \frac{z}{\eps^2}}, \label{eq:uepsplus}
\end{align}
where 
$\big(a_e(z),a_o(z)\big)$ is a random, Markovian process defined at $z \ge 0$, with initial condition  $a_t(0) = a_{t}^{(0)}$  given in \eqref{eq:Amplo} 
for $t \in \{e,o\}$, and with the infinitesimal generator  in Theorem  \ref{prop:1}.
To describe the results, it suffices to state 
from there that
\begin{equation}
|a_e(z)|^2 + |a_o(z)|^2 = \big(|a_{e}^{(0)}|^2 + |a_{o}^{(0)}|^2\big) \exp( -\Lambda z),
\label{eq:detpower}
\end{equation}
with probability one, and that 
\begin{equation}
\EE \big[ a_o(z) \overline{a_e(z)} \big] = a_{o}^{(0)} \overline{a_{e}^{(0)}}
\exp \big( -(\Gamma +\Lambda )z \big),
\label{eq:crossmoment}
\end{equation}
with positive $\Lambda$ and $\Gamma$ defined by 
\begin{align}
\Lambda &= 
\frac{k^4(n^2-1)^2D^2}{2\beta\big(\frac{2}{\eta} + D\big)} {\rm cos}^2\Big(\frac{\xi D}{2}\Big)
\int_0^{k^2} \hspace{-0.05in}\frac{d \gamma}{\pi \eta_\gamma \sqrt{\gamma}}
\Bigg[\frac{ 2\frac{\xi_\gamma^2}{\eta_\gamma^2}+ \sin^2(\xi_\gamma D) \big(1-\frac{\xi_\gamma^2}{\eta_\gamma^2}\big)}
{ 4\frac{\xi_\gamma^2}{\eta_\gamma^2}+\sin^2(\xi_\gamma D) \big(1-\frac{\xi_\gamma^2}{\eta_\gamma^2}\big)^2}\Bigg]
\widehat{\cal R} (\beta-\sqrt{\gamma}), \label{eq:defLambda}\\
\Gamma &= \frac{k^4(n^2-1)^2D^2}{ \beta^2\big(\frac{2}{\eta} + D\big)^2} {\rm cos}^4 \Big(\frac{\xi D}{2}\Big)
\widehat{\cal R}(0), \label{eq:defGamma}
\end{align}
in terms of 
\begin{equation}
\xi_\gamma = \sqrt{k^2 n^2 -\gamma}, \qquad \eta_\gamma  = \sqrt{k^2-\gamma},
\label{eq:defxietagamma}
\end{equation}
and the power spectral density $\widehat {\cal R} \ge 0$, the Fourier
transform of the covariance \eqref{eq:covar}.

We have therefore from \eqref{eq:uepsplus} and
\eqref{eq:detpower} that the total power of the guided waves decays
exponentially, at the rate $\Lambda$,
\begin{equation}
|u_+(z)|^2 + |u_-(z)|^2 = 2 \big(|a_e(z)|^2 +
|a_o(z)|^2\big)^2 = 2
\big(|a_{e}^{(0)}|^2 + |a_{o}^{(0)}|^2\big) \exp( -\Lambda z).
\label{eq:powerdecay}
\end{equation}
This decay models the transfer of power from the guided
modes to the radiation modes, induced by scattering at the random interfaces
\eqref{eq:Interfaces}. 

The imbalance of power
between the two waveguides is quantified by
\begin{align}
{\cal P}(z) &=\frac{|u_+(z)|^2-|u_-(z)|^2}{|u_+(z)|^2+|u_-(z)|^2 } = \frac{2 {\rm Re} \big\{  a_e(z)\overline{a_o(z)} \exp \big[i ( \beta_{e}- \beta_{o})  \frac{z}{\eps^2}\big] \big\}}{|a_e(z)|^2+|a_o(z)|^2 },
\end{align}
and its expectation is
\begin{align}
\label{eq:powertransf}
\EE\Big[\mathcal{P} (z)\Big] =\frac{2 {\rm Re}\big\{ a_{e}^{(0)}
  \overline{a_{o}^{(0)}} \exp\big( i ( \beta_{e}- \beta_{o})
  \frac{z}{\eps^2}\big)\big\} }{|a_{e}^{(0)}|^2+|a_{o}^{(0)}|^2 } \exp( -\Gamma z) .
\end{align}
To explain what this gives, consider the source  excitation \eqref{eq:assumea0}
with $a_{e}^{(0)}=a_{o}^{(0)}$, so that the initial wavefield is 
supported in the waveguide centered at $x = d/2 + D/2$. Then, equation \eqref{eq:powertransf} becomes
\begin{align}
\EE\Big[{\cal P} (z)\Big] =
\cos \Big[( \beta_{e}-  \beta_{o}) \frac{z}{\eps^2}\Big]
\exp( -\Gamma z) , 
\label{eq:imb1}
\end{align}
and it describes the competition between the deterministic
and random coupling of the waveguides.
The cosine in \eqref{eq:imb1} models the deterministic coupling which
induces periodic oscillations of the power,  as in section \ref{sect:resHomog}.
The random coupling is modeled by the exponential decay in $z$,
at the rate $\Gamma$.  It shows that as the range increases, the power tends to become equally distributed 
among the two waveguides. This decay is present in \eqref{eq:powertransf} as well, so the equipartition of power 
at large $z$ is  independent of the initial condition generated by the source.

\subsubsection{Weak coupling}
\label{sect:resWeak}
When the separation distance $d$ between the waveguides satisfies \eqref{eq:larged}, the wavefield 
$p(z/\eps^2,x)$ has the same expression as in \eqref{eq:pmoder}, \eqref{eq:uepsplus}, 
but the  random processes $(u_+(z),u_-(z))$ have different statistics. 

The total power of the guided waves is still given by \eqref{eq:powerdecay}, and decays at the 
same rate $\Lambda$ defined in \eqref{eq:defLambda}. However, the expectation of the imbalance of power between the two waveguides
satisfies the damped harmonic oscillator equation
\begin{align}
\label{eq:harmosc}
\big[\partial_z^2 +2\Gamma \partial_z +(2\theta \beta')^2 \big] \EE[\mathcal{P} (z)]  =0 ,
\end{align}
with $\beta'$ defined in \eqref{eq:Betaprime} and 
\begin{equation}
\label{eq:weak3}
\theta = \eps^{-2} \exp(-\eta d).
\end{equation}
Based on the value of $\theta$, which is independent of $\eps$ by assumption \eqref{eq:larged}, we distinguish three regimes, 
which we describe for the source excitation \eqref{eq:assumea0}, with $u_-^{(\eps)}(0) = 0$ and therefore $\EE[\mathcal{P} (0)] = 1$:

\vspace{0.05in} 
1. When $2 \theta \beta' < \Gamma$, the solution of \eqref{eq:harmosc} is 
\begin{align}
\nonumber
\EE[\mathcal{P} (z)] =& e^{-\Gamma z} 
\Big[ 
\cosh \big( \sqrt{\Gamma^2- (2\theta \beta' )^2} z\big)+
\frac{\Gamma}{\sqrt{\Gamma^2- (2\theta \beta' )^2}}
\sinh \big( \sqrt{\Gamma^2- (2\theta \beta' )^2} z\big) \Big].
\label{eq:weak4}
\end{align}
This tends to $1$ as $ \theta \to 0$, meaning that when the two
waveguides are very far apart, there is no transfer of power between them.
This is just as in the ideal (deterministic)  waveguide system. However, unlike in  ideal waveguides, 
the power is transferred from the guided modes to the radiation ones, as 
described by the exponential decay in  \eqref{eq:powerdecay}.

For a finite $\theta$, we have $\EE[\mathcal{P} (z)] \to 0$ as 
$z \to \infty$, so the random coupling distributes the power evenly among the two waveguides.

\vspace{0.05in} 
2. In the critical case $2 \theta \beta' = \Gamma$ we have 
\begin{equation}
\EE[\mathcal{P} (z)] = ( 1+ \Gamma z)\exp(-\Gamma z) .
\label{eq:weak5}
\end{equation}
As in the previous  case,   the random coupling equidistributes the power among the two waveguides, in the limit $z \to \infty$.

\vspace{0.05in} 
3. When $2 \theta \beta' > \Gamma$,  the solution of \eqref{eq:harmosc} is 
 \begin{align}
\nonumber
\EE[\mathcal{P} (z)]= e^{-\Gamma z}\Big[
\cos \big( \sqrt{ (2\theta \beta' )^2-\Gamma^2} z\big)+ 
\frac{\Gamma}{\sqrt{ (2\theta \beta' )^2-\Gamma^2}}
\sin \big( \sqrt{ (2\theta \beta' )^2-\Gamma^2} z\big) 
\Big]. \label{eq:weak6}
\end{align}
It displays periodic oscillations induced by the deterministic coupling of the 
waveguides, but these oscillations are damped due to the random coupling. 
In particular, if $2 \theta \beta' \gg \Gamma$, we get 
\begin{equation}
\EE[\mathcal{P} (z) ] \approx e^{-\Gamma z} \cos (2 \theta \beta' z),
\label{eq:weak7}
\end{equation}
in agreement with \eqref{eq:imb1}.

\begin{figure}
\centerline{ \includegraphics[width=6.5cm]{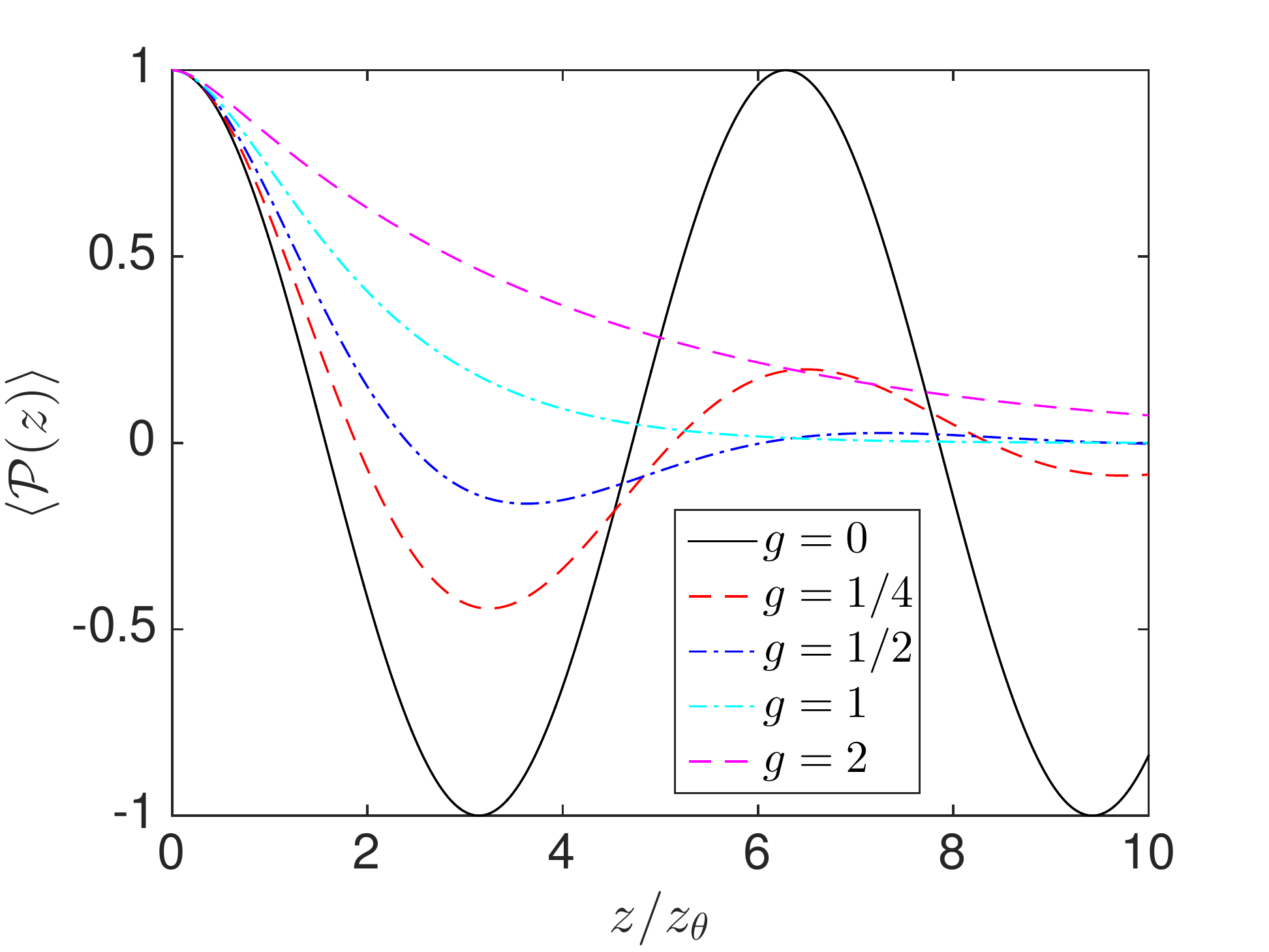} }
\caption{Imbalance ratio $\left<\mathcal{P} (z) \right> := \EE[{\cal P}(z)]$ as a function of $z/z_{\theta}$, 
where $z_{\theta} = 1/(2\theta \beta')$. We illustrate the result for different values of $g = \Gamma z_{\theta}$, given in the legend.
Note how the effective coupling coefficient $\Gamma$ reduces the deterministic and periodic transfer of power and causes the 
imbalance ratio  to tend to $0$.}
\label{fig:imbalance}
\end{figure}

\hspace{0.05in}
We plot in Figure \ref{fig:imbalance} the imbalance ratio $\EE[\mathcal{P} (z)]$  as a function of $z$ normalized 
by the deterministic coupling distance\footnote{At the distance $\pi z_\theta$ the power is fully transferred 
from one step-index waveguide to the other one, when there are no random perturbations.}  $z_{\theta} = 1/(2\theta \beta')$. 
From this plot and from equation \eqref{eq:powerdecay} we conclude that 
wave scattering at the random interfaces \eqref{eq:Interfaces}  has two net effects: 
\begin{enumerate}
\item  It induces  a self-averaging loss
(or leakage) of total power, due to the coupling of the guided  modes with the radiation modes. 
\item  It causes a blurring of the periodic (deterministic) power transfer from one waveguide to the other. \end{enumerate}
The blurring effect at item 2 is the main practical result of the paper. 
It depends on the effective parameter $\Gamma$ and it is significant as soon as $\Gamma z_\theta$ becomes of order one.
Therefore, the deterministic transfer of power is very sensitive to the random fluctuations of the 
interfaces \eqref{eq:Interfaces}.

\subsubsection{Very weak coupling}
\label{sect:resVWeak}
When the distance $d$ between the waveguides is so large that \eqref{eq:huged} holds, the wavefield has the same expression 
as \eqref{eq:pmoder}, but the wave amplitudes have different statistics. They model the  only coupling in this regime, 
between the guided and radiation modes, which generates effective wave power leakage.
Explicitly, we show in section \ref{sect:randomTwo} that the wave amplitudes converge in probability, as $\eps \to 0$,  
to the deterministic function $(|u_+(z)|^2,|u_-(z)|^2)$ satisfying 
\begin{equation}
\label{eq:vweakres}
|u_+(z)|^2 = |u_+^{(0)}(0)|^2\exp ( - \Lambda z ),\quad 
|u_-(z)|^2 = |u_-^{(0)}(0)|^2 \exp ( - \Lambda z ) .
\end{equation}
Although deterministic, this function is not as  in the ideal  waveguide system, where $(|u_+^{(0)}(z)|^2,|u_-^{(0)}(z)|^2)$ is constant in $z$
(because $\beta' e^{-\eta d} z \ll 1$ in (\ref{eq:up0}-\ref{eq:um0}) in this regime). Instead, it
decays exponentially at rate $\Lambda$, due to 
the power leakage. The  imbalance of power between the waveguides is constant
\begin{equation}
\label{eq:imb3}
{\cal P}(z) =  \frac{|u_+(z)|^2 - |u_-(z)|^2}{ |u_+(z)|^2+|u_-(z)|^2}  = \frac{|u_+^{(0)}(0)|^2 - |u_-^{(0)}(0)|^2}{ |u_+^{(0)}(0)|^2+|u_-^{(0)}(0)|^2} ,
\end{equation}
so in the case of the source excitation \eqref{eq:assumea0}, ${\cal P}(z) \simeq {\cal P}(0) = 1.$

\section{Analysis in ideal waveguides}
\label{sect:homog}

The analysis in this section is classical and follows the lines of \cite{magnanini}. It derives the results stated in section \ref{sect:resHomog}
by expanding the wave field  on a complete set of eigenmodes. The proof of the completeness of this set 
is the most delicate part and it is carried out in  \cite{magnanini} by the Levitan-Levinson method \cite[Chapter 9]{coddington}.

Recall the Helmholtz operator \eqref{eq:HelmIdeal} in the transverse coordinate, with index of refraction ${\rm n}^{(0)}(x)$ given in 
\eqref{eq:no}, and note that it is self-adjoint with respect to the scalar product associated to the $L^2$-norm,
\begin{align}
\label{scalarproduct}
(\phi_1,\phi_2) &:= \int_\RR  \overline{\phi_1(x)} {\phi_2(x)}dx .
\end{align}
Its spectrum is
$
(-\infty,  k^2) \cup \big\{\big(\beta_{t,j}^2\big)_{1 \le j \le N_t},  ~~ t \in \{e,o\}\big\}  ,
$
where $\beta_{t,j}$ are called the guided mode wavenumbers.
They are positive and satisfy  the order relation
\begin{equation}
k^2<\beta^2_{t,N_t} < \cdots <\beta_{t,1}^2 <k^2n^2, \end{equation}
where the index $t\in \{e,o\}$ stands for the even and odd eigenfunctions in the transverse coordinate $x$.

We describe next the eigenfunctions for the discrete spectrum  and the 
improper eigenfunctions for the continuum spectrum, and explain that they form a complete 
set. We use them to decompose the wavefield into guided, radiation and evanescent modes
with amplitudes determined by the source. 

\subsection{Discrete spectrum}
\label{sect:discrSpect}

There are two sets of discrete eigenvalues and eigenfunctions: the first one associated with the even modes in $x$ and 
the second one associated with the odd modes.

The $j$-th even eigenfunction $\phi_{e,j}(x)$ for the eigenvalue $\beta_{e,j}^2$ is  defined by 
\begin{align}
\frac{\phi_{e,j}(x)}{A_{e,j}} = 
\left\{
\begin{array}{ll}
 \exp(-\eta_{e,j} \frac{d}{2}) \sin (\xi_{e,j} D) \Big(1 +\frac{\xi_{e,j}^2}{\eta_{e,j}^2}\Big)\cosh (\eta_{e,j} x), \quad  & x \in \big[0, \frac{d}{2}\big] , \\
 \frac{\xi_{e,j}}{\eta_{e,j}}\cos \big[\xi_{e,j}(x-\frac{d}{2}-D)\big]-\sin\big[\xi_{e,j}(x-\frac{d}{2}-D)\big], &  x \in \big[\frac{d}{2},\frac{d}{2}+D\big] ,\\
\frac{\xi_{e,j}}{\eta_{e,j}} \exp\big[- \eta_{e,j} (x-\frac{d}{2}-D)\big], &  x \in \big[ \frac{d}{2}+D , \infty) ,
\end{array}
\right.
\label{eq:phie}
\end{align}
and  $\phi_{e,j}(-x) = \phi_{e,j}(x)$ for $x\geq 0$. Here 
\begin{align}
\xi_{e,j} = \sqrt{k^2n^2 -\beta_{e,j}^2}, \quad \eta_{e,j} = \sqrt{\beta_{e,j}^2-k^2} ,
\label{def:xiej}
\end{align}
and $A_{e,j}>0$ is the normalization constant 
\begin{align}
\nonumber
A_{e,j} =&\Bigg[
\frac{1}{2}\exp(-\eta_{e,j}d)\sin^2(\xi_{e,j} D) \Bigg(1+\frac{\xi_{e,j}^2}{\eta_{e,j}^2}\Bigg)^2
\Bigg( \frac{\sinh(\eta_{e,j} d)}{\eta_{e,j}}+d\Bigg)
\\
&+\Bigg(\frac{\xi_{e,j}^2}{\eta_{e,j}^2}-1\Bigg)\frac{\sin(2\xi_{e,j}D)}{2\xi_{e,j}}+\Bigg(\frac{\xi_{e,j}^2}{\eta_{e,j}^2}+1\Bigg)D
+2 \frac{\sin^2(\xi_{e,j}D)}{\eta_{e,j}}
+\frac{\xi_{e,j}^2}{\eta_{e,j}^3}\Bigg]^{-1/2} 
\end{align}
calculated so that $\phi_{e,j}$ has unit $L^2$-norm. Moreover,  $\beta_{e,j} \in (k , n k)$ satisfies the 
dispersion relation 
\begin{align}
\label{eq:dispersione}
\left(1+\frac{\xi_{e,j}^2}{\eta_{e,j}^2}\right)
e^{-\eta_{e,j} d} =  \left[1-\frac{\xi_{e,j}}{\eta_{e,j}}\tan\Big(\frac{\xi_{e,j} D}{2}\Big) \right]
\left[1+\frac{\xi_{e,j}}{\eta_{e,j}}{\rm cotan}\Big(\frac{\xi_{e,j} D}{2}\Big) \right] ,
\end{align}
with $\xi_{e,j},\eta_{e,j}$ given by (\ref{def:xiej}), which ensures the continuity of \eqref{eq:phie} and its derivative. 
The number of solutions $\beta_{e,j} \in (k , n k)$ of \eqref{eq:dispersione} is denoted by $N_e$ and when $d$ is large, it depends on the 
value of $k D \sqrt{n^2-1}$.

The $j$-th odd eigenfunction for the eigenvalue $\beta_{o,j}^2$ is defined similarly, 
\begin{align}
\frac{\phi_{o,j}(x)}{A_{o,j}} = 
\left\{
\begin{array}{ll}
\exp(-\eta_{o,j} \frac{d}{2}) \sin (\xi_{o,j} D) \left(1 +\frac{\xi_{o,j}^2}{\eta_{o,j}^2}\right)\sinh (\eta_{o,j} x), \quad & x \in \big[0,\frac{d}{2}\big] , \\
\frac{\xi_{o,j}}{\eta_{o,j}}\cos \big[\xi_{o,j}(x-\frac{d}{2}-D)\big] -\sin\big[\xi_{o,j}(x-\frac{d}{2}-D)\big], & x \in \big[\frac{d}{2}, \frac{d}{2}+D \big],\\
 \frac{\xi_{o,j}}{\eta_{o,j}} \exp\big[- \eta_{o,j} (x-\frac{d}{2}-D)\big], &  x \in \big[\frac{d}{2}+D , \infty \big),
\end{array}
\right.
\label{eq:phio}
\end{align}
with $\phi_{o,j}(-x) = - \phi_{o,j}(x)$ for $x\geq 0$ and 
\begin{align}
\label{eq:xio}
\xi_{o,j} = \sqrt{k^2n^2 -\beta_{o,j}^2}, \quad \eta_{o,j} =\sqrt{\beta_{o,j}^2-k^2}.
\end{align}
The normalization constant $A_{o,j}>0$ is given by
\begin{align}
\nonumber
A_{o,j} &=\Bigg[
\frac{1}{2} \exp(-\eta_{o,j}d) \sin^2(\xi_{o,j} D) \Bigg(1+\frac{\xi_{o,j}^2}{\eta_{o,j}^2}\Bigg)^2
\Bigg( \frac{\sinh(\eta_{o,j} d)}{\eta_{o,j}}-d\Bigg)
\\
&+\Bigg(\frac{\xi_{o,j}^2}{\eta_{o,j}^2}-1\Bigg)\frac{\sin(2\xi_{o,j}D)}{2\xi_{o,j}}+\Bigg(\frac{\xi_{o,j}^2}{\eta_{o,j}^2}+1\Bigg)D
+2\frac{\sin^2(\xi_{o,j}D)}{\eta_{o,j}}
+ \frac{\xi_{o,j}^2}{\eta_{o,j}^3} \Bigg]^{-1/2}
\label{eq:amplo}
\end{align}
so that $\phi_{o,j}$ has unit $L^2$ norm and $\beta_{o,j} \in ( k , n k)$  satisfies the dispersion relation
\begin{align}
\label{eq:dispersiono}
-\Bigg(1+\frac{\xi_{o,j}^2}{\eta_{o,j}^2}\Bigg)
e^{-\eta_{o,j} d} = \Bigg[1-\frac{\xi_{o,j}}{\eta_{o,j}}\tan\Big(\frac{\xi_{o,j} D}{2}\Big) \Bigg]
\Bigg[1+\frac{\xi_{o,j}}{\eta_{o,j}}{\rm cotan}\Big(\frac{\xi_{o,j} D}{2}\Big) \Bigg],
\end{align}
with $\xi_{o,j},\eta_{o,j}$ given by (\ref{eq:xio}), for  $j = 1, \ldots, N_o$, which ensures that \eqref{eq:phio} and its derivative are continuous.

\subsection{Continuous spectrum}
\label{sect:continSpect}%
For $\gamma \in (-\infty,k^2)$ there are two improper eigenfunctions, even and odd,  denoted by $\phi_{t,\gamma}(x)$, for $t \in \{e,o\}$. We write their expression below in terms of the parameters
\begin{align}
\label{def:xigamma}
\xi_\gamma = \sqrt{k^2n^2 -\gamma}, \quad \eta_\gamma = \sqrt{k^2-\gamma}.
\end{align}

The even eigenfunctions  satisfy $\phi_{e,\gamma}(-x) = \phi_{e,\gamma}(x)$ for $x \ge 0$, and are 
defined by 
\begin{align}
\frac{\phi_{e,\gamma}(x)}{A_{e,\gamma}} =
  \frac{\xi_\gamma}{\eta_\gamma} \cos (\eta_\gamma x), \quad  \mbox{for} ~ x \in \Big[0,\frac{d}{2}\Big],
\end{align}
and by
\begin{align}
\hspace{-0.1in}\frac{\phi_{e,\gamma}(x)}{A_{e,\gamma}} =  \frac{\xi_\gamma}{\eta_\gamma} \cos \Big[\xi_\gamma \Big(x-\frac{d}{2} \Big)\Big]
\cos \Big(\frac{\eta_\gamma d}{2}\Big) - \sin \Big[\xi_\gamma \Big(x-\frac{d}{2} \Big)\Big]
\sin \Big(\frac{\eta_\gamma d}{2}\Big),\nonumber \\ \mbox{for} ~ x \in \Big[\frac{d}{2}, \frac{d}{2}+D\Big], 
 \end{align}
  and by 
 \begin{align}
 \nonumber
\hspace{-0.05in}\frac{\phi_{e,\gamma}(x)}{A_{e,\gamma}} = 
\cos \Big[ \eta_\gamma \Big(x - \frac{d}{2}-D\Big)\Big] \Bigg[\frac{\xi_\gamma}{\eta_\gamma} \cos (\xi_\gamma D)
\cos\Big( \frac{\eta_\gamma d}{2}\Big) &- \sin (\xi_\gamma D)
\sin\Big( \frac{\eta_\gamma d}{2}\Big)\Bigg]
\\ 
-\frac{\xi_\gamma}{\eta_\gamma} \sin \Big[ \eta_\gamma \Big(x - \frac{d}{2}-D\Big)\Big] \Bigg[\frac{\xi_\gamma}{\eta_\gamma} \sin (\xi_\gamma D)
\cos \Big(\frac{\eta_\gamma d}{2}\Big) &+ \cos (\xi_\gamma D)
\sin \Big(\frac{\eta_\gamma d}{2}\Big)\Bigg] \nonumber \\
&\mbox{for} ~ x \ge \frac{d}{2}+D, \qquad \end{align}
 with normalization constant $A_{e,\gamma}>0$ given by
\begin{align}
\nonumber
A_{e,\gamma} =& (2 \pi \eta_\gamma)^{-1/2}  
 \Bigg\{\Bigg[\frac{\xi_\gamma}{\eta_\gamma} \cos (\xi_\gamma D)
\cos\Big( \frac{\eta_\gamma d}{2}\Big) - \sin (\xi_\gamma D)
\sin\Big( \frac{\eta_\gamma d}{2}\Big)\Bigg]^2 \nonumber \\
&+ \frac{\xi_\gamma^2}{\eta_\gamma^2}
\Bigg[\frac{\xi_\gamma}{\eta_\gamma} \sin (\xi_\gamma D)
\cos \Big(\frac{\eta_\gamma d}{2}\Big) + \cos (\xi_\gamma D)
\sin \Big(\frac{\eta_\gamma d}{2}\Big)\Bigg]^2 \Bigg\}^{-1/2}.
  \label{def:Aegamma}
\end{align}

The odd eigenfunctions satisfy $\phi_{o,\gamma}(-x) = -\phi_{o,\gamma}(x)$ for $x \ge 0$, and are
defined  by 
\begin{align}
\frac{\phi_{o,\gamma}(x)}{A_{o,\gamma}}\phi_{o,\gamma} =
  \frac{\xi_\gamma}{\eta_\gamma} \sin (\eta_\gamma x), \quad  \mbox{for} ~ x \in \Big[0,\frac{d}{2}\Big],
\end{align}
and  by  
\begin{align}
\hspace{-0.1in}  \frac{\phi_{o,\gamma}(x)}{A_{o,\gamma}} =  \frac{\xi_\gamma}{\eta_\gamma} \cos 
\Big[\xi_\gamma \Big(x-\frac{d}{2} \Big)\Big]
\sin \Big(\frac{\eta_\gamma d}{2}\Big) + \sin \Big[\xi_\gamma \Big(x-\frac{d}{2} \Big)\Big]
\cos \Big(\frac{\eta_\gamma d}{2}\Big) \nonumber \\ \mbox{for} ~ x \in \Big[\frac{d}{2}, \frac{d}{2}+D\Big], 
 \end{align}
 and by 
 \begin{align}
 \nonumber
\frac{\phi_{o,\gamma}(x)}{A_{o,\gamma}} =  
\cos \Big[ \eta_\gamma \Big(x - \frac{d}{2}-D\Big)\Big] \Bigg[\frac{\xi_\gamma}{\eta_\gamma} \cos (\xi_\gamma D)
\sin\Big( \frac{\eta_\gamma d}{2}\Big) &+ \sin (\xi_\gamma D)
\cos\Big( \frac{\eta_\gamma d}{2}\Big)\Bigg]
\\ 
\hspace{-0.35in}-\frac{\xi_\gamma}{\eta_\gamma} \sin \Big[ \eta_\gamma \Big(x - \frac{d}{2}-D\Big)\Big] 
\Bigg[\frac{\xi_\gamma}{\eta_\gamma} \sin (\xi_\gamma D)
\sin \Big(\frac{\eta_\gamma d}{2}\Big) &- \cos (\xi_\gamma D)
\cos \Big(\frac{\eta_\gamma d}{2}\Big)\Bigg]\nonumber \\
&\mbox{for} ~ x \ge \frac{d}{2}+D, \qquad 
\end{align}
with normalization constant $A_{o,\gamma}>0$ given by
\begin{align}
\nonumber
A_{o,\gamma} =& (2 \pi \eta_\gamma)^{-1/2}  
 \Bigg\{\Bigg[\frac{\xi_\gamma}{\eta_\gamma} \cos (\xi_\gamma D)
\sin\Big( \frac{\eta_\gamma d}{2}\Big) + \sin (\xi_\gamma D)
\cos\Big( \frac{\eta_\gamma d}{2}\Big)\Bigg]^2 \nonumber \\
&+ \frac{\xi_\gamma^2}{\eta_\gamma^2}
\Bigg[\frac{\xi_\gamma}{\eta_\gamma} \sin (\xi_\gamma D)
\sin \Big(\frac{\eta_\gamma d}{2}\Big) - \cos (\xi_\gamma D)
\cos \Big(\frac{\eta_\gamma d}{2}\Big)\Bigg]^2 \Bigg\}^{-1/2}.
  \label{def:Aogamma}
\end{align}

\vspace{0.05in}
\begin{remark}
\label{rem.1} Note that the improper eigenfunctions $\phi_{t,\gamma}(x)$, for $t \in \{e,o\}$, are not in $L^2(\RR)$. 
However, we can define for any $\varphi \in L^2(\RR)$ 
\begin{align*}
\left( \phi_{t,\gamma},\varphi\right) = \lim_{M \to +\infty} \int_{0}^M \phi_{t,\gamma} (x) \varphi(x)dx+ \lim_{M \to +\infty} \int_{-M}^0 
\phi_{t,\gamma}(x)  \varphi(x) dx
,
\end{align*}
where the limit holds in $L^2(-\infty,k^2)$. 
The normalizing constants $A_{t,\gamma}$ given in \eqref{def:Aegamma} and \eqref{def:Aogamma} are such that, 
 for any test function $a \in L^2(-\infty,k^2)$,
$$
2 \int_0^M \Big| \int_{-\infty}^{k^2} \phi_{t,\gamma}(x)a(\gamma) d\gamma\Big|^2 dx 
\stackrel{M \to +\infty}{\longrightarrow} \int_{-\infty}^{k^2} |a(\gamma)|^2 d\gamma .
$$
They depend continuously on $\gamma$ and are  bounded on $(-\infty,k^2)$ with $A_{t,\gamma} \simeq \big(2 \pi \sqrt{|\gamma|}\big)^{-1/2}$ as 
$\gamma \to -\infty$. 
\end{remark}

\subsection{Completness}
\label{sect:complete}
The proof of completness is based on the Levitan-Levinson method \cite[Chapter 9]{coddington} and is the same as that 
in \cite[Section 2.3]{magnanini}. We state directly the result:

For any $\varphi\in L^2(\RR)$ we have the Parseval identity
\begin{align*}
\left( \varphi,\varphi \right)  = \sum_{t \in \{e,o\}}  \sum_{j=1}^{N_t} \big| \left( \phi_{t,j}, \varphi\right) \big|^2
+\sum_{t \in \{e,o\}}  \int_{-\infty}^{k^2} \big| \left( \phi_{t,\gamma},\varphi\right)\big|^2 d\gamma   ,
\end{align*}
and the map which assigns to $\varphi$ the coefficients of its spectral decomposition
\begin{align*}
\varphi \mapsto \Big( \left( \phi_{t,j} , \varphi \right)  ,j=1,\ldots,N_t;
 \left(\phi_{t,\gamma},\varphi\right), \gamma\in ({-\infty},{k^2}); t\in \{e,o\}\Big)
\end{align*}
is an isometry from $L^2(\RR)$ onto $\CC^{N_e+N_o} \times L^2(-\infty,k^2)^2$.
Moreover, there exists a resolution of the identity $\Pi$ of the operator \eqref{eq:HelmIdeal} such that: for any $\varphi\in L^2(\RR)$ and  any $-\infty\leq r\leq r'\leq +\infty$,
\begin{align*}
\nonumber
\Pi_{(r,r')}(\varphi)(x) = &\sum_{t \in \{e,o\}} \sum_{j=1}^{N_t} \left(\phi_{t,j},  \varphi   \right) \phi_{t,j}(x) {\bf 1}_{(r,r')}(\beta_{t,j}^2)
\\
&+ \sum_{t \in \{e,o\}}
\int_r^{{\rm min}(k^2, r')}\hspace{-0.1in} \left(  \phi_{t,\gamma},\varphi \right) \phi_{t,\gamma}(x) d\gamma \, {\bf 1}_{(r,+\infty)}(k^2) ,
\end{align*}
and for any $\varphi$ in the domain of ${\cal H}$,
\begin{align*}
\nonumber
\Pi_{(r,r')}({\cal H}\varphi)(x) = & \sum_{t \in \{e,o\}}\sum_{j=1}^{N_t} \beta_{t,j}^2 \left( \phi_{t,j},\varphi   \right) \phi_{t,j}(x) {\bf 1}_{(r,r')}(\beta_{t,j}^2) \\
&+ \sum_{t \in \{e,o\}}
\int_r^{{\rm min}(k^2,r')} \hspace{-0.1in}\gamma \left(  \phi_{t,\gamma}, \varphi\right) \phi_{t,\gamma}(x) d\gamma \, {\bf 1}_{(r,+\infty)}(k^2) .
\end{align*}

\subsection{Modal decomposition}
Using the results in the previous section, we can write the solution $p^{(0)}(z,x)$ of the Helmholtz equation \eqref{eq:homogHelm}   as 
a superposition of modes
\begin{align}
\nonumber
{p}^{(0)}(z,x) =&  \sum_{t \in \{e,o\}} \sum_{j=1}^{N_t} {p}_{t,j}^{(0)}(z) \phi_{t,j}(x)
+ \sum_{t \in \{e,o\}}\int_0^{k^2} {p}_{t,\gamma}^{(0)} (z) \phi_{t,\gamma}(x) d\gamma \\
&
+ \sum_{t \in \{e,o\}} \int_{-\infty}^0 {p}_{t,\gamma}^{(0)} (z) \phi_{t,\gamma}(x) d\gamma .
\label{eq:modalexpansion}
\end{align}
The first sum contains the guided modes, where $p_{t,j}^{(0)}(z)$ are one-dimensional waves that propagate along $z$, with wavenumber $\beta_{t,j}$,
\begin{equation}
\label{eq:modal1}
\partial_z^2 {p}_{t,j}^{(0)}(z) +\beta_{t,j}^2 {p}_{t,j}^{(0)}(z)=0, \qquad z \ne 0, \quad  \mbox{for} ~j=1,\ldots,N_t,\quad t\in \{e,o\}.
\end{equation}
The other two sums contain the radiation modes and  evanescent modes, where 
\begin{align}
\label{eq:modal2}
\partial_z^2 {p}_{t,\gamma}^{(0)}(z) +\gamma {p}_{t,\gamma}^{(0)}(z)&=0, \qquad z \ne 0, \quad  \mbox{for} ~\gamma \in (-\infty,k^2),\quad t\in \{e,o\} .
\end{align}
The radiation modes correspond to  $\gamma \in (0,k^2)$, so that ${p}_{t,\gamma}^{(0)}(z)$ are one-dimensional 
waves that propagate along $z$ at wavenumber $\sqrt{\gamma}$. The evanescent modes are decaying exponentially\footnote{
We neglect the case $\gamma = 0$ because there is no coupling of the modes in  ideal waveguides, 
and therefore it has a negligible contribution in \eqref{eq:modalexpansion}. Mode coupling occurs in randomly perturbed waveguides, so the case $\gamma = 0$ must be dealt with carefully, as explained in section \ref{sec:limitthm}.}
  in $z$ on the range scale $1/\sqrt{|\gamma|}$, for $\gamma < 0$.

The solutions of (\ref{eq:modal1}--\ref{eq:modal2})  that satisfy the radiation conditions are
\begin{align*}
p_{t,j}^{(0)}(z) &= \frac{a_{t,j}^{(0)}}{\sqrt{\beta_{t,j}}} e^{i \beta_{t,j} z} 1_{(0,\infty)}(z) + \frac{b_{t,j}^{(0)}}{\sqrt{\beta_{t,j}}} e^{-i \beta_{t,j} z} 1_{(-\infty,0)}(z), 
\quad &&j = 1, \ldots, N_t, 
\\
p_{t,\gamma}^{(0)}(z) &= \frac{a_{t,\gamma}^{(0)}}{\gamma^{1/4}} e^{i \sqrt{\gamma} z} 1_{(0,\infty)}(z) + \frac{b_{t,\gamma}^{(0)}}{\gamma^{1/4}} 
e^{-i \sqrt{\gamma} z} 1_{(-\infty,0)}(z), 
&&\gamma \in (0,k^2), 
 \\
p_{t,\gamma}^{(0)}(z) &= \frac{a_{t,\gamma}^{(0)}}{|\gamma|^{1/4}} e^{-\sqrt{|\gamma|} z} 1_{(0,\infty)}(z) + \frac{b_{t,\gamma}^{(0)}}{\gamma^{1/4}} 
e^{ \sqrt{|\gamma|} z} 1_{(-\infty,0)}(z), 
&& \gamma \in (\infty,0),
\end{align*}
for $t \in \{e,o\}$, with constant mode amplitudes determined by the source in (\ref{eq:fouriertransform}),
\begin{align}
\label{eq:sourceconda1}
{a}_{t,j}^{(0)} = - {b}_{t,j}^{(0)} =& \frac{\sqrt{\beta_{t,j}}}{2}  \left( \phi_{t,j} ,{f} \right)  ,\quad && \hspace{-0.6in}j=1,\ldots,N_t,\\
\label{eq:sourceconda2}
{a}_{t,\gamma}^{(0)} = - {b}_{t,\gamma}^{(0)} =& \frac{\gamma^{1/4}}{2} \left( \phi_{t,\gamma},{f}  \right), && \hspace{-0.6in} \gamma \in (0,k^2) ,\\
{a}_{t,\gamma}^{(0)} =- {b}_{t,\gamma}^{(0)} =& \frac{|\gamma|^{1/4}}{2} \left( \phi_{t,\gamma} ,{f}  \right), && \hspace{-0.6in}\gamma \in (-\infty,0) .
\label{eq:sourceconda3}
\end{align}
Substituting in \eqref{eq:modalexpansion}, we obtain the modal expansion of the wavefield $p^{(0)}(z,x)$.

\subsection{Single guided component waveguides}
\label{subsec:unmode}
We describe here the special situation considered in section \ref{sect:resHomog}, where the assumptions \eqref{eq:assumeD} and \eqref{eq:dlarge}
hold and there is a single guided mode  per step-index waveguide, as we now explain.

\vspace{0.05in}
\begin{lemma}
\label{lem.singlebeta}
Under the assumption $k D \sqrt{n^2-1} < \pi$, equation 
\begin{equation}
\label{eq:barbeta1}
\frac{\sqrt{ n^2 k^2-\beta^2}}{\sqrt{\beta^2-k^2}}\tan\big(\frac{\sqrt{ n^2 k^2-\beta^2} D}{2}\big)=1
\end{equation}
has a unique solution $\beta$ in $(k,n k)$, whereas
equation
\begin{equation}
\label{eq:barbeta2}
\frac{\sqrt{ n^2 k^2-\beta^2}}{\sqrt{\beta^2-k^2}}{\rm cotan}\big(\frac{\sqrt{ n^2 k^2-\beta^2} D}{2}\big) =-1
\end{equation}
has no solution in $(k,nk)$. Furthermore, when the distance $d$ between the waveguides is large enough, so that 
$\exp(-\eta d) \ll 1$, with $\eta = \sqrt{\beta^2-k^2}$, equations \eqref{eq:dispersione} and \eqref{eq:dispersiono} 
have a single solution, meaning that $N_e = N_o = 1$.
\end{lemma}

\vspace{0.1in}
The proof of the lemma is in appendix \ref{sect:ap0}. We let henceforth $\beta$ be the unique solution of \eqref{eq:barbeta1} in the interval $(k,nk)$ and 
recall the definition \eqref{eq:defxieta} of $\xi$ and $\eta$ in terms of $\beta$. Since $N_e = N_o = 1$, 
we simplify notation as $\beta_{t,1} \leadsto \beta_t,$ for  $t \in \{e,o\}$, and obtain from \eqref{eq:dispersione} and \eqref{eq:dispersiono}
that 
\begin{align}
\beta_{e} &= \beta + \beta' e^{-\eta d} +o\Big( e^{-\eta d }\Big), \label{eq:approxbetae} \qquad \beta_{o} = \beta - \beta' e^{-\eta d}+o\Big( e^{-\eta d }\Big),
\end{align}
with
\begin{align}
 \nonumber
 \beta' &= -\Bigg\{ \frac{\partial}{\partial \beta_s}\Bigg[ \frac{\sqrt{k^2n^2-\beta_s^2}}{\sqrt{\beta_s^2-k^2}} \tan\Bigg( \frac{\sqrt{k^2n^2-\beta_s^2} D}{2}\Bigg) \Bigg] \Bigg\}^{-1} \Bigg|_{\beta_s=\beta}= \frac{\eta}{\beta \big(1+\frac{\eta^2}{\xi^2}\big) \big(\frac{1}{\eta}+\frac{D}{2}\big)} ,
\end{align}
Similarly, we let $\phi_{t,1}(x) \leadsto \phi_t(x)$, for $t \in \{e,o\}$, 
and obtain from \eqref{eq:phie}--\eqref{eq:amplo} that 
\begin{align}
\phi_{e}(x) = \phi(x) +O\big( \exp(-\eta d )\big) , \quad 
\phi_{o}(x) = \phi(x) +O\big( \exp(-\eta d )\big) , \quad x\geq 0. \label{eq:phi_eo}
\end{align}
Substituting in \eqref{eq:modalexpansion} and obtaining from \cite[Section 3.2]{magnanini}  that the radiation and 
evanescent modes have an $O(z^{-2})$  contribution\footnote{Note that in [10] the decay is $O(z^{-1})$ because 
the source is $f(x) \delta(z)$. In our case the decay is $O(z^{-2})$,  because we have the $z$ derivative of the solution 
in [10], due to the source $f(x) \delta'(z)$.}, we obtain 
equations \eqref{eq:phom} and (\ref{eq:Po1}--\ref{eq:um0}) used in section  \ref{sect:resHomog} to describe the 
deterministic coupling between the ideal step-index waveguides. 
\section{Analysis in randomly perturbed waveguides}
\label{sect:random}
\begin{figure}[t]
\vspace{-0.6in}
\begin{picture}(0,0)%
\hspace{1.5in}\includegraphics[width = 0.45\textwidth]{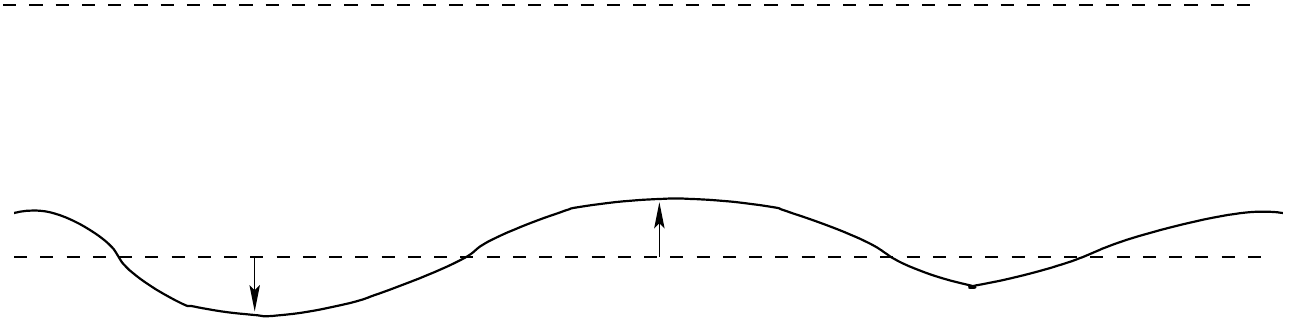}%
\end{picture}%
\setlength{\unitlength}{2763sp}%
\begingroup\makeatletter\ifx\SetFigFont\undefined%
\gdef\SetFigFont#1#2#3#4#5{%
  \reset@font\fontsize{#1}{#2pt}%
  \fontfamily{#3}\fontseries{#4}\fontshape{#5}%
  \selectfont}%
\fi\endgroup%
\begin{picture}(8819,2177)(1779,-3866)
\put(3000,-3700){\makebox(0,0)[lb]{\smash{{\SetFigFont{7}{8.4}{\familydefault}{\mddefault}{\updefault}{\color[rgb]{0,0,0}{\normalsize $x = -\frac{d}{2} - D$}}%
}}}}
\put(3000,-2900){\makebox(0,0)[lb]{\smash{{\SetFigFont{7}{8.4}{\familydefault}{\mddefault}{\updefault}{\color[rgb]{0,0,0}{\normalsize $x = -\frac{d}{2}$}}%
}}}}
\put(5000,-3500){\makebox(0,0)[lb]{\smash{{\SetFigFont{7}{8.4}{\familydefault}{\mddefault}{\updefault}{\color[rgb]{0,0,0}{\normalsize $z_1$}}%
}}}}
\put(6300,-3900){\makebox(0,0)[lb]{\smash{{\SetFigFont{7}{8.4}{\familydefault}{\mddefault}{\updefault}{\color[rgb]{0,0,0}{\normalsize $z_2$}}%
}}}}
\end{picture}%
\vspace{0.1in}
\caption{Illustration of fluctuations of the bottom interface at $x = -d/2-D$.   The ideal interfaces are drawn with the dotted lines, whereas the fluctuations are drawn with the full line. We identify two ranges  $z_1$ and $z_2$ so that 
$\nu_1(z_1) < 0$ and  $\nu_1(z_2) > 0$. }
\label{fig:interface}
\end{figure}
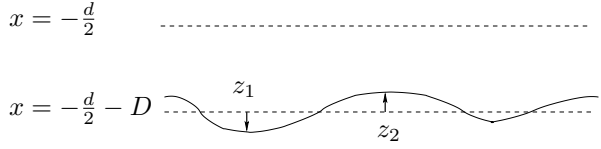

Let us introduce the notation 
\begin{equation}
\Delta k^2 = k^2(n^2-1)
\end{equation}
and  rewrite the index of refraction 
defined in (\ref{eq:modelpert1}--\ref{eq:Interfaces}) as 
\begin{align}
\label{eq:newindref}
\big({\rm n}^{(\eps)}(z,x)\big)^2 = \big({\rm n}^{(0)}(x)\big)^2+  (n^2-1)   V^{(\eps)}(z,x), \quad \quad V^{(\eps)}(z,x) =\sum_{q=1}^4 V_q^{(\eps)}(z,x) ,
\end{align}
with
\begin{align*}
V_1^{(\eps)}(z,x) = &  -  {\bf 1}_{(-d/2-D, -d/2-D+\eps D \nu_1(z))}(x) {\bf 1}_{(0,+\infty)}(\nu_1(z)) \\
&+ {\bf 1}_{(-d/2-D+\eps D \nu_1(z), -d/2-D)}(x) {\bf 1}_{(-\infty,0)}(\nu_1(z))   ,\\
V_2^{(\eps)}(z,x) = &   {\bf 1}_{(-d/2, -d/2+\eps D \nu_2(z))}(x) {\bf 1}_{(0,+\infty)}(\nu_2(z)) \\
&- {\bf 1}_{(-d/2+\eps D \nu_2(z), -d/2)}(x) {\bf 1}_{(-\infty,0)}(\nu_2(z))   ,\\
V_3^{(\eps)}(z,x) = &   - {\bf 1}_{(d/2, d/2+\eps D \nu_3(z))}(x) {\bf 1}_{(0,+\infty)}(\nu_3(z)) \\
&+ {\bf 1}_{(d/2+\eps D \nu_3(z), d/2)}(x) {\bf 1}_{(-\infty,0)}(\nu_3(z))    ,\\
V_4^{(\eps)}(z,x) = &   {\bf 1}_{(d/2+D, d/2+D+\eps D \nu_4(z))}(x) {\bf 1}_{(0,+\infty)}(\nu_4(z)) \\
&- {\bf 1}_{(d/2+D+\eps D\nu_4(z), d/2+D)}(x) {\bf 1}_{(-\infty,0)}(\nu_4(z))  .
\end{align*}
To explain this reformulation, we illustrate in  Figure \ref{fig:interface}  the fluctuations $\eps D \nu_1(z)$ of
 the interface at $x = - d/2 -D$. At 
range $z_1$ we have $\nu_1(z_1) < 0$, so   $(z_1,x)$ with $x \in (-d/2-D + \eps D \nu_1(z_1),-d/2-D)$ is in the waveguide and therefore, 
\[
\big({\rm n}^{(\eps)}(z_1,x)\big)^2  = \big({\rm n}^{(0)}(x)\big)^2 +(n^2-1)V_1^{(\eps)}(z_1,x) = 1 + (n^2-1) = n^2.
\]
At range $z_2$ we have $\nu_1(z_2) > 0$, so  $(z_2,x)$ with $x \in (-d/2-D ,-d/2-D+\eps D \nu_1(z))$ is outside the waveguide 
and 
\[
\big({\rm n}^{(\eps)}(z_2,x)\big)^2  = \big({\rm n}^{(0)}(x)\big)^2 -(n^2-1)V_1^{(\eps)}(z_2,x) = n^2 - (n^2-1) = 1.
\]
The same reasoning applies to all the interfaces.

Our goal in this section is to analyze the solution $p(z,x)$ of \eqref{eq:fouriertransform}, with index of refraction 
\eqref{eq:newindref} and radiation condition at infinity, and in particular, to derive the results stated in section \ref{sect:resRand}.

\subsection{Modal decomposition}
\label{sect:modecRand}

The completness result of section \ref{sect:complete}
allows  the expansion of  $p(z,x)$, $z$ by $z$, in terms of the eigenfunctions defined 
in sections \ref{sect:discrSpect}--\ref{sect:continSpect},
\begin{align}
\nonumber
{p}(z,x) =&  \sum_{t \in \{e,o\}} \sum_{j=1}^{N_t} {p}_{t,j}(z) \phi_{t,j}(x)
+ \sum_{t \in \{e,o\}}\int_0^{k^2} {p}_{t,\gamma}(z) \phi_{t,\gamma}(x) d\gamma \\
&
+ \sum_{t \in \{e,o\}} \int_{-\infty}^0 {p}_{t,\gamma} (z) \phi_{t,\gamma}(x) d\gamma .
\label{eq:modalexpansionRan}
\end{align}
Here $p_{t,j}(z)$ are complex-valued amplitudes of guided modes with wavenumber $\beta_{t,j}$ satisfying 
\eqref{eq:dispersione} or \eqref{eq:dispersiono}. They propagate along $z$ and satisfy the
following one-dimensional Helmholtz equations at $ z> 0$, 
\begin{align}\nonumber
\partial_z^2 {p}_{t,j}(z) +\beta_{t,j}^2 {p}_{t,j}(z)=& -\Delta k^2 \sum_{t' \in \{e,o\}}
\sum_{l'=1}^{N_{t'}} C_{t,j,t',l'}^{(\eps)}(z) {p}_{t',l'}(z)  \\&\hspace{-0.4in}-\Delta k^2 
 \sum_{t' \in \{e,o\}} \int_{-\infty}^{k^2} C_{t,j,t',\gamma'}^{(\eps)} (z) {p}_{t',\gamma'}(z) d\gamma', \quad  j=1,\ldots,N_t.
\label{eq:cma2a}
\end{align}
Similarly, 
$p_{t,\gamma}(z)$ are complex-valued amplitudes of  modes that are propagating for $\gamma \in (0,k^2)$ (radiation modes)
and  decaying for $\gamma < 0$ (evanescent modes), and satisfy the one-dimensional Helmholtz equations
\begin{align}
\nonumber
\partial_z^2 {p}_{t,\gamma}(z) +\gamma {p}_{t,\gamma}(z)=& -\Delta k^2 
 \sum_{t' \in \{e,o\}}\sum_{l'=1}^{N_{t'}} C_{t,\gamma,t' ,l'}^{(\eps)}(z) {p}_{t',l'}(z) \\
&
\hspace{-0.4in}-\Delta k^2 \sum_{t' \in \{e,o\}} \int_{-\infty}^{k^2} C_{t,\gamma,t',\gamma'}^{(\eps)} (z) {p}_{t',\gamma'}(z) d\gamma',  
\quad  \gamma \in (-\infty,k^2)  .
\label{eq:cma2b}
\end{align}
The source terms in these equations are due to the random fluctuations, which are supported at $z \in (0,L/\eps^2)$. These 
induce mode coupling, modeled by the  random coefficients 
\begin{align}
C_{t,j,t',l'}^{(\eps)}(z) =& \left(\phi_{t,j},\phi_{t',l'} V^{(\eps)}(z,\cdot) \right) ,\label{eq:defCoef1}\\
C_{t,j,t',\gamma'}^{(\eps)}(z) =& \left(\phi_{t,j},\phi_{t',\gamma'} V^{(\eps)}(z,\cdot) \right) ,\label{eq:defCoef2}\\
C_{t,\gamma,t', l'}^{(\eps)}(z) =&  \left(\phi_{t,\gamma},\phi_{t',l'} V^{(\eps)}(z,\cdot) \right) ,\label{eq:defCoef3}\\
C_{t,\gamma,t',\gamma'}^{(\eps)}(z) =& \left(\phi_{t,\gamma},\phi_{t',\gamma'} V^{(\eps)}(z,\cdot) \right) .\label{eq:defCoef4}
\end{align}
Recalling the definition \eqref{eq:newindref} of $V^{(\eps)}(z,x)$ and  using Taylor expansions of the eigenfunctions $\phi_{t,j}(x)$ and $\phi_{t,\gamma}(x)$ 
around $x = \pm d/2$ and $x = \pm (d/2 + D)$, we  obtain a power series (in $\eps$) expression
of these coefficients
\begin{align}
C_{t,j,t',l'}^{(\eps)} (z) =& \eps C_{t,j,t',l'}(z)  + \eps^2 c_{t,j,t',l'}(z)  +o(\eps^2) ,\\
\nonumber
C_{t,j,t',l'}(z) =& -  D \nu_1(z) [\phi_{t,j} \phi_{t',l'}] \Big(\hspace{-0.05in}-\frac{d}{2}-D\Big) + D \nu_2(z)[\phi_{t,j} \phi_{t',l'}] 
\Big(\hspace{-0.05in}-\frac{d}{2}\Big) \\
&-  D\nu_3(z)[\phi_{t,j} \phi_{t',l'}] \Big(\frac{d}{2}\Big) +  D\nu_4(z) [\phi_{t,j} \phi_{t',l'}] \Big(\frac{d}{2}+D\Big) 
\label{expres:Cttp} ,\\
\nonumber
c_{t,j,t',l'}(z) =& - \frac{D^2 \nu_1^2(z)}{2} \partial_x[ \phi_{t,j} \phi_{t',l'}] \Big(\hspace{-0.05in}-\frac{d}{2}-D\Big) +  \frac{D^2\nu_2^2(z)}{2} \partial_x[\phi_{t,j} \phi_{t',l'}] \Big(\hspace{-0.05in}-\frac{d}{2}\Big) \\
&-  \frac{D^2\nu_3^2(z)}{2} \partial_x[\phi_{t,j} \phi_{t',l'}] \Big(\frac{d}{2}\Big) + \frac{D^2\nu_4^2(z) }{2} \partial_x[ \phi_{t,j} \phi_{t',l'}] \Big(\frac{d}{2}+D\Big)   , \label{eq:csmall}
\end{align}
and similarly for  $C_{t,j,t',\gamma'}^{(\eps)}, C_{t,\gamma,t', l'}^{(\eps)} , C_{t,\gamma,t',\gamma'}^{(\eps)}$.

Since we consider propagation distances of the order of $\eps^{-2}$, we neglect from now on the $o(\eps^2)$ terms 
in these expansions, because they give no contribution in the limit $\eps \to 0$. To simplify the presentation, 
we will also not write the $\eps^2$ terms, even though they are not negligible. They contribute to the expression of   the infinitesimal generators
of the limit processes via terms that  are proportional to $\frac{\Delta k^2}{2 \beta_{t,j}} \EE [ c_{t,j,t,j}(z) ]$. 
Because $\partial_x (\phi_{t,j}^2)$ is odd, for $t \in \{e,o\}$, we obtain from \eqref{eq:csmall} that 
\begin{align*}
\EE [ c_{t,j,t,j}(z) ] &= \frac{D^2}{2}\EE[\nu_4^2(z) +\nu_1^2(z) ]  \partial_x \phi_{t,j}^2\Big(\frac{d}{2}+D\Big) 
-
\frac{D^2}{2} \EE[\nu_3^2(z) + \nu_2^2(z) ]\partial_x \phi_{t,j}^2\Big(\frac{d}{2}\Big)    \\
&=
D^2{\cal R}(0)\Big[  \partial_x \phi_{t,j}^2\Big(\frac{d}{2}+D\Big)   
- \partial_x \phi_{t,j}^2\Big(\frac{d}{2}\Big)\Big] .
\end{align*}
Moreover, in the case of two weakly coupled single guided component waveguides,
\begin{equation}
\label{eq:DetPhase}
\EE [ c_{t,1,t,1}(z) ] =- 2 D^2 {\cal R}(0) \xi \sin(\xi D) \Big(\frac{2}{\eta} + D\Big)^{-1}.
\end{equation}
We will incorporate these  terms in the statements of the results
without giving additional  details of their derivation. They
correspond to effective deterministic phases of the mode amplitudes, in the limit $\eps\to0$, as explained in Remark \ref{rem:3}.

\subsection{Forward and backward going waves}
The guided modes can be decomposed further in forward and backward guided modes. 
This decomposition is basically the  method of variation of parameters for the perturbed Helmholtz equations 
(\ref{eq:cma2a}--\ref{eq:cma2b}), where we define the complex valued amplitudes 
\begin{equation}
\label{eq:amplitudes}
\{a_{t,j}(z), \, b_{t,j}(z), ~ j = 1, \ldots, N_t\} ~~ \mbox{and} ~~ \{a_{t,\gamma}(z), \, b_{t,\gamma}(z), ~ \gamma \in (0,k^2)\},
\end{equation}
for $t \in \{e,o\}$, such that 
\begin{align}
{p}_{t,j}(z) =& \frac{1}{\sqrt{\beta_{t,j}}}\Big( {a}_{t,j} (z) e^{i\beta_{t,j} z} +{b}_{t,j}(z) e^{- i\beta_{t,j} z} \Big), \nonumber \\
\partial_z  {p}_{t,j}(z) =& i\sqrt{\beta_{t,j}}\Big( {a}_{t,j} (z) e^{i\beta_{t,j}z} -{b}_{t,j}(z) e^{- i\beta_{t,j} z} \Big),\quad j=1,\ldots,N_t,
\label{eq:guidedFB}
\end{align}
and 
\begin{align}
{p}_{t,\gamma}(z) =& \frac{1}{\gamma^{1/4}}\Big( {a}_{t,\gamma} (z) e^{i\sqrt{\gamma} z} 
+{b}_{t,\gamma}(z) e^{- i\sqrt{\gamma}z} \Big), \nonumber \\
\partial_z  {p}_{t,\gamma}(z) =& i\gamma^{1/4}\Big( {a}_{t,\gamma} (z) e^{i\sqrt{\gamma} z} 
-{b}_{t,\gamma}(z) e^{- i\sqrt{\gamma}z} \Big),\quad \gamma\in (0,k^2) .
\label{eq:radFB}
\end{align}
Substituting (\ref{eq:guidedFB}--\ref{eq:radFB}) in (\ref{eq:cma2a}--\ref{eq:cma2b}),  we obtain that 
the amplitudes \eqref{eq:amplitudes}
 satisfy the following first-order system of stochastic differential equations 
\begin{align}
\nonumber
\partial_z {a}_{t,j}(z) =& \frac{i \eps\Delta k^2}{2} \hspace{-0.1in}\sum_{t' \in \{e,o\}}
\sum_{l'=1}^{N_{t'}} \frac{C_{t,j,t',l'}(z)}{\sqrt{\beta_{t',l'}\beta_{t,j}}} \Big[ {a}_{t',l'}(z) e^{i (\beta_{t',l'}-\beta_{t,j})z}
+ {b}_{t',l'}(z) e^{i (-\beta_{t',l'}-\beta_{t,j})z}
\Big]\\
\nonumber
&\hspace{-0.4in}+  \frac{i \eps\Delta k^2}{2} \hspace{-0.1in} \sum_{t' \in \{e,o\}} 
\int_0^{k^2}  \frac{C_{t,j,t',\gamma'}(z)}{\sqrt[4]{\gamma'}\sqrt{\beta_{t,j}}} \Big[ {a}_{t',\gamma'}(z) e^{i (\sqrt{\gamma'}-\beta_{t,j})z}
+ {b}_{t',\gamma'}(z) e^{i (-\sqrt{\gamma'}-\beta_{t,j})z}
\Big]d\gamma'\\
&\hspace{-0.4in}+
 \frac{i \eps\Delta k^2}{2} \hspace{-0.1in}\sum_{t' \in \{e,o\}}
\int_{-\infty}^0 \frac{C_{t,j,t',\gamma'}(z)}{\sqrt{\beta_{t,j}}}  {p}_{t',\gamma'}(z) e^{-i \beta_{t,j}z}
d\gamma' ,
\label{eq:evol1a}
\end{align}
\begin{align}
\nonumber
\partial_z {a}_{t,\gamma}(z) =& \frac{i \eps\Delta k^2}{2} \hspace{-0.1in}\sum_{t' \in \{e,o\}}
\sum_{l'=1}^{N_{t'}} \frac{C_{t,\gamma,t',l'}(z)}{\sqrt[4]{\gamma} \sqrt{\beta_{t',l'}}} \Big[ {a}_{t',l'}(z) e^{i (\beta_{t',l'}-\sqrt{\gamma})z}
+ {b}_{t',l'}(z) e^{i (-\beta_{t',l'}-\sqrt{\gamma})z}
\Big]\\
\nonumber
&\hspace{-0.4in}+  \frac{i \eps\Delta k^2}{2} \hspace{-0.1in}\sum_{t' \in \{e,o\}} 
\int_0^{k^2}  \frac{C_{t,\gamma,t',\gamma'}(z)}{\sqrt[4]{{\gamma'}{\gamma}}} \Big[ {a}_{t',\gamma'}(z) e^{i (\sqrt{\gamma'}-\sqrt{\gamma})z}
+ {b}_{t',\gamma'}(z) e^{i (-\sqrt{\gamma'}-\sqrt{\gamma})z}
\Big]d\gamma'\\
&\hspace{-0.4in}+
 \frac{i \eps\Delta k^2}{2} \hspace{-0.1in} \sum_{t' \in \{e,o\}}
\int_{-\infty}^0 \frac{C_{t,\gamma,t',\gamma'}(z)}{\sqrt[4]{\gamma}}  {p}_{t',\gamma'}(z) e^{-i \sqrt{\gamma}z}
d\gamma' ,
\label{eq:evol1b}
\end{align}
\begin{align}
\nonumber
\partial_z {b}_{t,j}(z) =& -\frac{i \eps\Delta k^2}{2} \hspace{-0.1in}\sum_{t' \in \{e,o\}}
\sum_{l'=1}^{N_{t'}} \frac{C_{t,j,t',l'}(z)}{\sqrt{\beta_{t',l'}\beta_{t,j}}} \Big[ {a}_{t',l'}(z) e^{i (\beta_{t',l'}+\beta_{t,j})z}
+ {b}_{t',l'}(z) e^{i (-\beta_{t',l'}+\beta_{t,j})z} \Big]\\
\nonumber
&\hspace{-0.4in}-  \frac{i \eps\Delta k^2}{2}\hspace{-0.1in} \sum_{t' \in \{e,o\}} 
\int_0^{k^2}  \frac{C_{t,j,t',\gamma'}(z)}{\sqrt[4]{\gamma'}\sqrt{\beta_{t,j}}} \Big[ {a}_{t',\gamma'}(z) e^{i (\sqrt{\gamma'}+\beta_{t,j})z}
+ {b}_{t',\gamma'}(z) e^{i (-\sqrt{\gamma'}+\beta_{t,j})z}
\Big] d\gamma'\\
&\hspace{-0.4in}-
 \frac{i \eps\Delta k^2}{2} \hspace{-0.1in}\sum_{t' \in \{e,o\}}
\int_{-\infty}^0 \frac{C_{t,j,t',\gamma'}(z)}{\sqrt{\beta_{t,j}}}  {p}_{t',\gamma'}(z) e^{i \beta_{t,j}z}
d\gamma' ,
\label{eq:evol1c}
\end{align}
\begin{align}
\nonumber
\partial_z {b}_{t,\gamma}(z) =&- \frac{i \eps\Delta k^2}{2} \hspace{-0.1in}\sum_{t' \in \{e,o\}}
\sum_{l'=1}^{N_{t'}} \frac{C_{t,\gamma,t',l'}(z)}{\sqrt[4]{\gamma}\sqrt{\beta_{t',l'}\}}} \Big[ {a}_{t',l'}(z) e^{i (\beta_{t',l'}+\sqrt{\gamma})z}
+ {b}_{t',l'}(z) e^{i (-\beta_{t',l'}+\sqrt{\gamma})z}
\Big]\\
\nonumber
&\hspace{-0.4in}-  \frac{i\eps \Delta k^2}{2}\hspace{-0.1in} \sum_{t' \in \{e,o\}} 
\int_0^{k^2}  \frac{C_{t,\gamma,t',\gamma'}(z)}{\sqrt[4]{{\gamma'}{\gamma}}} \Big[ {a}_{t',\gamma'}(z) e^{i (\sqrt{\gamma'}+\sqrt{\gamma})z}
+ {b}_{t',\gamma'}(z) e^{i (-\sqrt{\gamma'}+\sqrt{\gamma})z}
\Big]d\gamma'\\
&\hspace{-0.4in}-
 \frac{i \eps\Delta k^2}{2} \hspace{-0.1in}\sum_{t' \in \{e,o\}}
\int_{-\infty}^0 \frac{C_{t,\gamma,t',\gamma'}(z)}{\sqrt[4]{\gamma}}  {p}_{t',\gamma'}(z) e^{i \sqrt{\gamma}z}
d\gamma'.
\label{eq:evol1d}
\end{align}
The right-hand sides in these equations are supported at $z \in (0,L/\eps^2)$ and model the mode coupling induced by scattering at the random interfaces \eqref{eq:Interfaces}. This coupling causes the randomization of the amplitudes \eqref{eq:amplitudes} and 
therefore of the wavefield. 

\subsection{Role of the evanescent modes}
Equations (\ref{eq:evol1a}--\ref{eq:evol1d}) show that the forward and backward going amplitudes of both the guided and radiation modes are coupled to each other and to the evanescent modes $p_{t,\gamma}(z)$, for $\gamma \in (-\infty,0)$ and $t \in \{e,o\}$. These mode amplitudes satisfy 
\begin{align}
\label{eq:evan1}
\partial_z^2 {p}_{t,\gamma}(z) +\gamma {p}_{t,\gamma}(z)= 
- \eps \big[ g_{t,\gamma}(z) + g^{\rm ev}_{t,\gamma}(z) \big], \quad z \ne 0, 
\end{align}
where $g_{t,\gamma}(z)$ and $g^{\rm ev}_{t,\gamma}(z)$ are supported at $z \in (0,L/\eps^2)$ and are defined by 
\begin{align}
\nonumber
g_{t,\gamma}(z) =& \Delta k^2\hspace{-0.1in} \sum_{t' \in \{e,o\}}
\sum_{l'=1}^{N_{t'}} \frac{C_{t,\gamma,t',l'}(z)}{\sqrt{\beta_{t',l'}}} \Big[ {a}_{t',l'}(z) e^{i \beta_{t',l'} z}
+ {b}_{t',l'}(z) e^{- i\beta_{t',l'} z}
\Big]\\
&+   \Delta k^2  \hspace{-0.1in}\sum_{t' \in \{e,o\}} 
\int_0^{k^2}  \frac{C_{t,\gamma,t',\gamma'}(z)}{\sqrt[4]{\gamma'}} \Big[ {a}_{t',\gamma'}(z) e^{i  \sqrt{\gamma'} z}
+ {b}_{t',\gamma'}(z) e^{- i \sqrt{\gamma'} z}
\Big]d\gamma',  \label{eq:g_tg}\\
g^{\rm ev}_{t,\gamma}(z) =&
 \Delta k^2 \hspace{-0.1in} \sum_{t' \in \{e,o\}}
\int_{-\infty}^0 C_{t,\gamma,t',\gamma'}(z) {p}_{t',\gamma'}  (z)
d\gamma' , \qquad z \in \Big(0,{L}/{\eps^2}\Big).
\label{eq:defgtga}
\end{align}

We now explain that the solution of \eqref{eq:evan1}  can be expressed in terms of the  amplitudes \eqref{eq:amplitudes} when $\eps$ is small enough. This result has already been obtained in similar frameworks \cite{gomez11} and it allows 
us to obtain a closed system of equations for the guided and radiation mode amplitudes.

Let us invert the  operator $\partial_z^2 + \gamma$ in \eqref{eq:evan1}  using the  Green's function that satisfies 
the radiation condition (i.e., it  decays away from the source). This gives  
\begin{align}
{p}_{t,\gamma} (z) =& 
\frac{\eps}{2\sqrt{|\gamma|}}
\int_{0}^{{L}/\eps^2} \big[g_{t,\gamma}(z') + g^{\rm ev}_{t,\gamma}(z')\big]  e^{-\sqrt{|\gamma|} |z-z'|} dz'
+ \frac{{a}_{t,\gamma}^{(0)}}{ {|\gamma|}^{1/4}} e^{-\sqrt{|\gamma|} z}  ,
\label{eq:GFint}
\end{align}
at $ z > 0$, for $\gamma<0$ and $t\in \{e,o\}$.
The first term in this expression comes from the random mode coupling induced by scattering at the random 
interfaces \eqref{eq:Interfaces}.
The second term  comes from the source.
Since we will consider propagation distances $z$ of the order of $\eps^{-2}$, we can neglect  this second term.

Using \eqref{eq:defgtga} and  the notation ${\itbf p}^{\rm ev}_t =({p}_{t,\gamma}(z))_{z>0,\gamma <0}$ for $t \in \{e,o\}$, we obtain that 
\begin{align}
\big( I- \eps {\cal J}_t \big) {\itbf p}^{\rm ev}_t  = \eps {\itbf G}_t ,
\end{align}
where  $I $ is the identity operator,  ${\cal J}_t$ is the integral operator
\begin{align}
[{\cal J}_t{\itbf p}^{\rm ev}_t ]_\gamma(z) =& \frac{\Delta k^2}{2\sqrt{|\gamma|}}\hspace{-0.03in}  \sum_{t' \in \{e,o\}}
\int_{0}^{L/\eps^2} dz' \hspace{-0.05in}
\int_{-\infty}^0 d \gamma' C_{t,\gamma,t',\gamma'}(z') {p}_{t',\gamma'} (z') 
 e^{-\sqrt{|\gamma|} |z-z'|}
 \end{align}
coming from the $g_{t,\gamma}^{\rm ev}$ term in \eqref{eq:GFint}, and the right-hand side  is  
\begin{equation}
[ {\itbf G}_t  ]_\gamma(z)=   \frac{\Delta k^2}{2\sqrt{|\gamma|}}
\int_{0}^{{L}/\eps^2}  
{g}_{t,\gamma}(z') e^{-\sqrt{|\gamma|} |z-z'|} dz'  , \label{eq:defBigG}
\end{equation}
for $z>0,\gamma < 0$.
The operator ${\cal J}_t$ is bounded from ${\cal C}^0( (0,\infty), L^1(-\infty,0))$ to itself and 
as in \cite{gomez11}, it can be shown that 
\[
\lim_{\eps \to 0}\PP( I- \eps {\cal J}_t \mbox{ is invertible})=1,
\]
and therefore 
\begin{align}
 {\itbf p}^{\rm ev}_t = \eps {\itbf G}_t 
+O(\eps^2).
\end{align}
Furthermore, equations (\ref{eq:evol1a}--\ref{eq:evol1d}) give that ${a}_{t',l'}(z')$, 
${b}_{t',l'}(z')$, ${a}_{t',\gamma'}(z')$, ${b}_{t',\gamma'}(z')$ are equal to 
${a}_{t',l'}(z)$, 
${b}_{t',l'}(z)$, ${a}_{t',\gamma'}(z)$, ${b}_{t',\gamma'}(z)$ up to terms of order $\eps$, as long as $|z-z'|=O(1)$, where the 
exponential in \eqref{eq:defBigG} contributes.

Gathering the results we obtain 
\begin{align}
& {p}_{t,\gamma}(z) =  
\frac{\eps \Delta k^2}{2\sqrt{|\gamma|}}   \sum_{t' \in \{e,o\}}
\int_{0}^{L/\eps^2} 
\sum_{l'=1}^{N_{t'}} \Bigg\{\frac{C_{t,\gamma,t',l'}(z')}{\sqrt{\beta_{t',l'}}} \Big[ {a}_{t',l'}(z) e^{i \beta_{t',l'} z'}
+ {b}_{t',l'} (z)e^{- i\beta_{t',l'} z'}
\Big]\nonumber \\
&+\int_0^{k^2}  \frac{C_{t,\gamma,t',\gamma'}(z')}{\sqrt[4]{\gamma'}} \Big[ {a}_{t',\gamma'}(z) e^{i  \sqrt{\gamma'} z'}
+ {b}_{t',\gamma'}(z) e^{- i \sqrt{\gamma'} z'}
\Big] d\gamma'  \Bigg\}
e^{-\sqrt{|\gamma|} |z-z'|} dz'   + O(\eps^2), 
\end{align}
for $z > 0, \gamma < 0$ and $t \in \{e, o\}$. 
Substituting this expression into (\ref{eq:evol1a}-\ref{eq:evol1d}) we get a closed system
for the guided and radiation mode amplitudes \eqref{eq:amplitudes}.

\subsection{Long range scaling}
The resulting stochastic system for the  guided and radiation mode amplitudes has a right-hand side
which is proportional to $\eps$. Therefore, the amplitudes are not affected by scattering at the random interfaces
\eqref{eq:Interfaces}  until the waves travel at long enough range,  dependent of $\eps$.  Note that the 
$O(\eps)$ coupling terms in (\ref{eq:evol1a}--\ref{eq:evol1d}) have zero expectation, so we need a central limit theorem type scaling, 
with range $z \leadsto  {z}/{\eps^2}, $ in order to observe net scattering effects. 

We rename the mode amplitudes in this long range scaling as 
\begin{align}
&
{a}_{t,j}^{(\eps)}(z) = {a}_{t,j}\Big(\frac{z}{\eps^2}\Big),\quad
{b}_{t,j}^{(\eps)}(z) = {b}_{t,j}\Big(\frac{z}{\eps^2}\Big)\, \quad j=1,\ldots,N_t ,\nonumber \\
&{a}_{t,\gamma}^{(\eps)}(z) = {a}_{t,\gamma}\Big(\frac{z}{\eps^2}\Big),~~\,
{b}_{t,\gamma}^{(\eps)}(z) = {b}_{t,\gamma}\Big(\frac{z}{\eps^2}\Big), \quad \gamma\in (0,k^2),
\label{eq:rescAmplitudes}
\end{align}
for $t\in \{e,o\}$.
Recalling that the random fluctuations are supported in the range interval $(0,L/\eps^2)$ and using the radiation 
conditions (i.e., the waves are outgoing at  $z > L$), we obtain
\begin{align}
& {a}_{t,j}^{(\eps)}(z=0) = {a}_{t,j}^{(0)} ,\quad 
{b}_{t,j}^{(\eps)}(z=L) = 0, \quad j=1,\ldots,N_t, 
\label{eq:boundCond1} \\
&{a}_{t,\gamma}^{(\eps)}(z=0) = {a}_{t,\gamma}^{(0)},  ~~\,
{b}_{t,\gamma}^{(\eps)}(z=L) = 0 , \quad \gamma\in (0,k^2),
 \label{eq:boundCond2}\end{align}
where ${a}_{t,j}^{(0)} $ and ${a}_{t,\gamma}^{(0)}$ are given in
(\ref{eq:sourceconda1}--\ref{eq:sourceconda2}), for $t\in \{e,o\}$ .  These boundary conditions and the closed  
system of stochastic differential equations described at the end of the previous section define the random amplitudes
\eqref{eq:rescAmplitudes}. 

\subsection{Forward scattering approximation}
We now introduce  the forward scattering approximation, where 
the coupling between forward and backward going modes can be
neglected. It follows from the following facts, under the assumption that the power spectral density $\widehat {\cal R}(\kappa)$ has compact support or fast decay
so that 
\begin{equation}
\widehat {\cal R}(\kappa) = \int_{-\infty}^\infty {\cal R}(z) e^{i \kappa z} dz \approx 0, \quad \mbox{if}~ ~ |\kappa| \geq  k.
\label{eq:assumeForward}
\end{equation}
\vspace{0.05in} 
1.  In the limit $\eps \to 0$, the coupling between the 
forward and backward guided modes depends on the coefficients
\begin{equation}
\int_0^\infty \EE\big[ C_{t,j,t',l'}(0)C_{t,j,t',l'}(z) \big] \cos\big[(\beta_{t',l'}+\beta_{t,j})z\big] dz  ,
\label{eq:coefftype1}
\end{equation}
while the coupling among the forward-guided modes depends
on
\begin{equation}
\int_0^\infty \EE\big[ C_{t,j,t',l'}(0)C_{t,j,t',l'}(z) \big] \cos\big[(\beta_{t',l'}-\beta_{t,j})z\big] dz  .
\label{eq:coefftype2}
\end{equation}
Recalling definitions \eqref{eq:newindref} and (\ref{eq:defCoef1}--\ref{eq:defCoef4}) and that $\beta_{t,j}
\in (k,n k)$, we conclude that the coefficients \eqref{eq:coefftype1}, which are proportional to 
$\widehat {\cal R}(\beta_{t',l'}+\beta_{t,j})$, are negligible under the assumption \eqref{eq:assumeForward}.
Therefore,  we can neglect the coupling between the  forward and backward guided modes. 
Nevertheless,   the forward guided modes are coupled among themselves, because the coefficients \eqref{eq:coefftype2},
which are proportional to  $\widehat {\cal R}(\beta_{t',l'}-\beta_{t,j})$, are  not negligible.

\vspace{0.05in} 2. The coupling between forward guided and
backward radiation modes depends in the limit $\eps \to 0$ on  the coefficients
\begin{equation}
\int_0^\infty \EE\big[ C_{t,j,t',\gamma'}(0)C_{t,j,t',\gamma'}(z) \big] \cos\big[(\sqrt{\gamma'}+\beta_{t,j})z\big] dz  ,
\end{equation}
for $\gamma' \in (0,k^2)$. These are proportional to $\widehat {\cal R}(\beta_{t,j}+\sqrt{\gamma'})$ and since $\beta_{t,j} \in (k,nk)$,
this type of coupling is negligible under the assumption \eqref{eq:assumeForward}.

\vspace{0.05in}
3. The coupling  between any radiation modes
is negligible in our regime $\eps \to 0$, as we will see in the following.
This is because these modes are essentially spatially supported outside the region 
where the medium is fluctuating.

\vspace{0.05in}
We suppose henceforth that assumption \eqref{eq:assumeForward} holds.
Because there is no coupling between the forward and backward going modes, and the backward mode amplitudes 
satisfy the homogeneous boundary conditions (\ref{eq:boundCond1}--\ref{eq:boundCond2}) at $z = L$, we can 
set 
\[
b_{t,j}^{(\eps)}(z) \approx 0, \quad b_{t,\gamma}^{(\eps)}(z) \approx 0, \qquad j = 1, \ldots, N_t, ~~ \gamma \in (0,k^2), ~~ 
t \in \{e,o\}.
\]
The expression \eqref{eq:modalexpansionRan} of the wavefield simplifies to 
\begin{align}
{p}\Big(\frac{z}{\eps^2},x\Big) =
\hspace{-0.05in}\sum_{t \in \{e,o\}} \Bigg[ \sum_{j=1}^{N_t} \frac{{a}_{t,j}^{(\eps)} (z) }{\sqrt{\beta_{t,j}}}  e^{i\beta_{t,j} \frac{z}{\eps^2}} \phi_{t,j}(x)+
\int_0^{k^2}  \frac{{a}_{t,\gamma}^{(\eps)} (z)}{\gamma^{1/4}}  e^{i\sqrt{\gamma} \frac{z}{\eps^2}} \phi_{t,\gamma}(x)d\gamma 
\Bigg]
+o(1) ,
\label{eq:expresspeps}
\end{align}
and the forward guided mode amplitudes ${a}_{t,j}^{(\eps)} $ and ${a}_{t,\gamma}^{(\eps)}$ satisfy 
\begin{align}
\nonumber
\partial_z {a}_{t,j}^{(\eps)}(z) =& \frac{i\Delta k^2}{2 \eps} \sum_{t' \in \{e,o\}}
\sum_{l'=1}^{N_{t'}} \frac{C_{t,j,t',l'}(\frac{z}{\eps^2})}{\sqrt{\beta_{t',l'}\beta_{t,j}}}  {a}_{t',l'}^{(\eps)}(z) e^{i (\beta_{t',l'}-\beta_{t,j}) \frac{z}{\eps^2}} \\
\nonumber
&+  \frac{i\Delta k^2}{2 \eps} \sum_{t' \in \{e,o\}} 
\int_0^{k^2}  \frac{C_{t,j,t',\gamma'}(\frac{z}{\eps^2})}{\sqrt[4]{\gamma'}\sqrt{\beta_{t,j}}}   {a}_{t',\gamma'}^{(\eps)}(z) e^{i (\sqrt{\gamma'}-\beta_{t,j}) \frac{z}{\eps^2}}
d\gamma'\\
&+
 \frac{i \Delta k^2}{2} \sum_{t' \in \{e,o\}}
\int_{-\infty}^0 \frac{C_{t,j,t',\gamma'}(\frac{z}{\eps^2})}{\sqrt{\beta_{t,j}}}  {q}_{t',\gamma'}^{(\eps)} (z)e^{-i \beta_{t,j} \frac{z}{\eps^2}}
d\gamma' , \label{eq:ForwEq1}
\end{align}
and 
\begin{align}
\nonumber
\partial_z {a}_{t,\gamma}^{(\eps)}(z) =& \frac{i \Delta k^2}{2 \eps} \sum_{t' \in \{e,o\}}
\sum_{l'=1}^{N_{t'}} \frac{C_{t,\gamma,t',l'}( \frac{z}{\eps^2} )}{\sqrt{\beta_{t',l'}}\sqrt[4]{\gamma}}   {a}_{t',l'}^{(\eps)}(z) e^{i (\beta_{t',l'}-\sqrt{\gamma}) \frac{z}{\eps^2} } \\
\nonumber
&+  \frac{i \Delta k^2}{2 \eps} \sum_{t' \in \{e,o\}} 
\int_0^{k^2}  \frac{C_{t,\gamma,t',\gamma'}(\frac{z}{\eps^2})}{\sqrt[4]{{\gamma'}{\gamma}}}  {a}_{t',\gamma'}^{(\eps)}(z) e^{i (\sqrt{\gamma'}-\sqrt{\gamma})\frac{z}{\eps^2}}
d\gamma'\\
&+
 \frac{i \Delta k^2}{2} \sum_{t' \in \{e,o\}}
\int_{-\infty}^0 \frac{C_{t,\gamma,t',\gamma'}(\frac{z}{\eps^2})}{\sqrt[4]{\gamma}}  {q}_{t',\gamma'}^{(\eps)}(z) e^{-i \sqrt{\gamma} \frac{z}{\eps^2}}
d\gamma' ,\label{eq:ForwEq2}
\end{align}
where 
\begin{align}
\nonumber
{q}_{t,\gamma}^{(\eps)} (z)=&  
\frac{\Delta k^2}{2\sqrt{|\gamma|}}   \sum_{t' \in \{e,o\}}
\int_{0}^{{L}/\eps^2}  \Bigg[ 
\sum_{l'=1}^{N_{t'}} \frac{C_{t,\gamma,t',l'}(z')}{\sqrt{\beta_{t',l'}}}   {a}_{t',l'}^{(\eps)}(z) e^{i \beta_{t',l'} z'}
 \\
&+\int_0^{k^2}  \frac{C_{t,\gamma,t',\gamma'}(z')}{\sqrt[4]{\gamma'}} {a}_{t',\gamma'}^{(\eps)}(z) e^{i  \sqrt{\gamma'} z'}
 d\gamma'  \Bigg]
e^{-\sqrt{|\gamma|}\big| \frac{z}{\eps^2}-z'\big|} dz' .
\end{align}

\subsection{Diffusion limit theorem}
\label{sec:limitthm}%
The next theorem gives the $\eps \to 0$ limit of the random forward-going mode amplitudes, 
the solutions of (\ref{eq:ForwEq1}--\ref{eq:ForwEq2}) with the initial conditions given in 
(\ref{eq:boundCond1}--\ref{eq:boundCond2}).
\ 
\begin{theorem}
\label{prop:1}
Suppose that the wavenumbers $\beta_{t,j}$  are distinct, for $1 \le j \le N_t$ and $t \in \{e,o\}$. Then, the  process 
$
\big( ({a}_{t,j}^{(\eps)}(z) )_{j=1}^{N_t}, ({a}_{t,\gamma}^{(\eps)}(z) )_{\gamma \in (0,k^2)},  t\in \{e,o\}\big)
$
converges in distribution in ${\cal C}^0\big([0,{L}], \CC^{N}  \times L^2(0,k^2)^2 \big)$, where $\CC^{N}  \times L^2(0,k^2)^2$
is equipped with the weak topology,  to the Markov process  
$
\big( ({a}_{t,j}(z) )_{j=1}^{N_t} , ({a}_{t,\gamma} (z))_{\gamma \in (0,k^2)},  t\in \{e,o\} \big).
$
The infinitesimal generator of this limit process is 
$
{\cal L}= {\cal L}^1+{\cal L}^2+{\cal L}^3,
$
where $\{{\cal L}^j\}_{1 \le j \le 3}$ are the differential operators:
\begin{align}
\nonumber {\cal L}^1 = & \frac{1}{2}\sum_{t,t' \in \{e,o\}}
\sum_{j=1}^{N_t} \sum_{l'=1}^{N_{t'}}
\Gamma^{c}_{t,j,t',l'} \big( 
{a}_{t,j} \overline{{a}_{t,j}} \partial_{{a}_{t',l'}} \partial_{\overline{{a}_{t',l'}}}
+
{a}_{t',l'} \overline{{a}_{t',l'}} \partial_{{a}_{t,j}} \partial_{\overline{{a}_{t,j}}} \\
\nonumber &\quad \quad
-
{a}_{t,j}  {a}_{t',l'} \partial_{{a}_{t,j}} \partial_{ {a}_{t',l'}}
-
\overline{{a}_{t,j}}  \overline{{a}_{t',l'}} \partial_{\overline{{a}_{t,j}}} \partial_{ \overline{{a}_{t',l'}}}\big)
{\bf 1}_{(t,j)\neq(t',l')}
\\
\nonumber &+
 \frac{1}{2}\sum_{t,t' \in \{e,o\}}
\sum_{j=1}^{N_t} \sum_{l'=1}^{N_{t'}}
\Gamma^1_{t,j,t',l'} \big( 
{a}_{t,j} \overline{{a}_{t',l'}} \partial_{{a}_{t,j}} \partial_{\overline{{a}_{t',l'}}}
+
\overline{{a}_{t,j}} {a}_{t',l'} \partial_{\overline{{a}_{t,j}}} \partial_{{a}_{t',l'}}\\
\nonumber &\quad \quad-
{a}_{t,j}  {a}_{t',l'} \partial_{{a}_{t,j}} \partial_{ {a}_{t',l'}}
-
\overline{{a}_{t,j}}  \overline{{a}_{t',l'}} \partial_{\overline{{a}_{t,j}}} \partial_{ \overline{{a}_{t',l'}}}\big)
\\
\nonumber &+ \frac{1}{2}\sum_{t\in \{e,o\}}
\sum_{j=1}^{N_t} \big( \Gamma^{c}_{t,j,t,j} - \Gamma^1_{t,j,t,j}\big)
\big( {a}_{t,j} \partial_{{a}_{t,j}} + \overline{{a}_{t,j}} \partial_{\overline{{a}_{t,j}}}
\big)\\
&
+\frac{i}{2}
\sum_{t\in \{e,o\}}
\sum_{j=1}^{N_t}  \Gamma^{s}_{t,j,t,j}  
\big( {a}_{t,j} \partial_{{a}_{t,j}} - \overline{{a}_{t,j}} \partial_{\overline{{a}_{t,j}}}
\big)  , \label{eq:defL1}
\\
{\cal L}^2
=&
-\frac{1}{2} \sum_{t\in \{e,o\}}
\sum_{j=1}^{N_t} ( \Lambda^{c}_{t,j}  +i\Lambda^{s}_{t,j} ) {a}_{t,j} \partial_{{a}_{t,j}} 
+
( \Lambda^{c}_{t,j}  - i\Lambda^{s}_{t,j} )  \overline{{a}_{t,j}} \partial_{\overline{{a}_{t,j}}}  ,
 \label{eq:defL2}\\{\cal L}^3
=&
i   \sum_{t\in \{e,o\}}
\sum_{j=1}^{N_t} \big( \kappa_{t,j} + \kappa_{t,j}^{\rm ev} \big)   \big( {a}_{t,j} \partial_{{a}_{t,j}} 
- \overline{{a}_{t,j}} \partial_{\overline{{a}_{t,j}}} \big) . \label{eq:defL3}
\end{align}
In these definitions we use the  classical complex derivative:
\[
 \zeta=\zeta_r+i\zeta_i~ \leadsto ~\partial_\zeta=\frac{1}{2}(\partial_{\zeta_r}-i \partial_{\zeta_i}), \quad 
\partial_{\overline{\zeta}} =\frac{1}{2}(\partial_{\zeta_r} +i \partial_{\zeta_i}),
\]
and the coefficients of the operators \eqref{eq:defL1}--\eqref{eq:defL3} are defined for 
indexes $j = 1, \ldots, N_t,$ $l' = 1, \ldots, N_{t'}$ and $t,t' \in \{e,o\}$, as follows: 

\vspace{0.05in} \noindent For all $(t,j)\neq(t',l')$:
\begin{align*}
\Gamma^{c}_{t,j,t',l'} =&
\frac{\Delta k^4}{2 \beta_{t,j}\beta_{t',l'}}
\int_0^\infty \EE\big[ C_{t,j,t',l'}(0) C_{t,j,t',l'}(z) \big] \cos \big[ (\beta_{t',l'}-\beta_{t,j})z\big] dz ,
\\
\Gamma^{s}_{t,j,t',l'} =&
\frac{\Delta k^4}{2 \beta_{t,j}\beta_{t',l'}}
\int_0^\infty \EE\big[ C_{t,j,t',l'}(0) C_{t,j,t',l'}(z) \big] \sin \big[(\beta_{t',l'}-\beta_{t,j})z\big] dz .
\end{align*}
For all $(t,j),(t',l')$:
\begin{align*}
\Gamma^1_{t,j,t',l'} =&
\frac{\Delta k^4}{2 \beta_{t,j}\beta_{t',l'}}
\int_0^\infty \EE\big[ C_{t,j,t,j}(0) C_{t',l',t',l'}(z) \big]   dz.
\end{align*}
For all $(t,j)$:
\begin{align*}
\Gamma^{c}_{t,j,t,j} =& -\hspace{-0.05in}\sum_{l =1,l\neq j}^{N_t} \Gamma^{c}_{t,j,t,l} - \sum_{t' \neq t} \sum_{l' =1}^{N_{t'}}
\Gamma^{c}_{t,j,t',l'} ,\\
\Gamma^{s}_{t,j,t,j} =& -\hspace{-0.05in}\sum_{l =1,l\neq j}^{N_t} \Gamma^{s}_{t,j,t,l} - \sum_{t' \neq t} \sum_{l' =1}^{N_{t'}}
\Gamma^{s}_{t,j,t',l'} ,\\
\Lambda_{t,j}^{c} =& \hspace{-0.05in}\sum_{t'\in \{e,o\}} \int_0^{k^2}  \hspace{-0.05in} \frac{\Delta k^4}{2\sqrt{\gamma}\beta_{t,j}}
\int_0^\infty \EE \big[ C_{t,j,t',\gamma}(0) C_{t,j,t',\gamma}(z) \big] 
\cos \big[ (\sqrt{\gamma} - \beta_{t,j})z\big] dz d\gamma,
\label{def:Lambdac} \\
\Lambda_{t,j}^{s} =&\hspace{-0.05in} \sum_{t'\in \{e,o\}}  \int_0^{k^2} \hspace{-0.05in}\frac{\Delta k^4}{2\sqrt{\gamma}\beta_{t,j}}
\int_0^\infty \EE \big[ C_{t,j,t',\gamma}(0) C_{t,j,t',\gamma}(z) \big] 
\sin \big[ (\sqrt{\gamma}-\beta_{t,j})z\big] dz d\gamma, \\
\kappa_{t,j}^{\rm ev} = &  \hspace{-0.05in} \sum_{t'\in \{e,o\}}  \int_{-\infty}^0  \frac{\Delta k^4}{2\sqrt{|\gamma|}\beta_{t,j}}
\int_0^\infty \EE \big[ C_{t,j,t',\gamma}(0) C_{t,j,t',\gamma}(z) \big] \cos ( \beta_{t,j} z) e^{-\sqrt{|\gamma|}z}dz d\gamma , \\
\kappa_{t,j} = &  \frac{\Delta k^2 D^2}{2\beta_{t,j}}
{\cal R}(0) \Bigg[ \partial_x( \phi_{t,j}^2)\Big(\frac{d}{2}+D\Big)   -
 \partial_x ( \phi_{t,j}^2)\Big(\frac{d}{2}\Big) \Bigg].
\end{align*} 
\end{theorem}

The proof of this diffusion limit theorem is based on a martingale approach using the perturbed test function method. 
It is an extension of the diffusion approximation theorem in \cite{papa74,fouque} that is carried out in \cite{gomez11} and \cite[Sections 3, 4.1]{gomez09b}.
Here we reproduce the result but we simplify its statement. The rigorous statement \cite[Theorem 4.1]{gomez09b} has two steps:
The first step deals with the convergence of truncated processes, obtained using the  resolution of the identity $\Pi_{\RR \setminus (-\gamma_\star,\gamma_\star)}$ defined in section \eqref{sect:complete}.  This  removes the vicinity of $\gamma = 0$, in order to take 
the diffusion limit $\eps \to 0$. The second step deals with the convergence of the truncated processes themselves, in the limit $\gamma_\star \to 0$. 

\begin{remark}
\label{rem:2}
Note that: 

\vspace{0.05in} 1. The convergence result in Theorem \ref{prop:1} holds in the weak topology. 

\vspace{0.05in} 2. The infinitesimal generator ${\cal L}$ does not involve the derivatives $\partial_{a_{t,\gamma}}$. Therefore 
$\Big( ({a}_{t,j}^{(\eps)} (z))_{j=1}^{N_t}, t\in \{e,o\}  \Big)$
converges in distribution in ${\cal C}^0([0,{L}], \CC^{N}  )$ to
the Markov process  $\Big( ({a}_{t,j}(z) )_{j=1}^{N_t} ,t\in \{e,o\}   \Big)$
with  generator ${\cal L}$. Of course,  the weak and strong topologies are the same  in $\CC^N$.\\

\vspace{0.05in} 3. If the generator ${\cal L}$ is applied to a test function that depends only on the mode powers $\Big( (P_{t,j} = |{a}_{t,j}|^2 )_{j=1}^{N_t} ,t\in \{e,o\}   \Big)$,
then the result is a function that depends only on $\Big( (P_{t,j} )_{j=1}^{N_t} ,t\in \{e,o\}   \Big)$. Thus, 
the mode powers $\Big( (P_{t,j}(z) )_{j=1}^{N_t} ,t\in \{e,o\}   \Big)$ define  a Markov process, with infinitesimal generator
\begin{align}
\nonumber
{\cal L}_{{\itbf P}} =  &\sum_{t,t' \in \{e,o\}}
\sum_{j=1}^{N_t} \sum_{l'=1}^{N_{t'}}
\Gamma^{c}_{t,j,t',l'} \Big( P_{t',l'} P_{t,j} \big( \partial_{P_{t_j}} -\partial_{P_{t',l'}}\big)
\partial_{P_{t,j}} + \big( P_{t',l'}-P_{t,j}\big) \partial_{P_{t,j}} \Big)\\
&
-  \sum_{t\in \{e,o\}}
\sum_{j=1}^{N_t}  \Lambda^{c}_{t,j} P_{t,j}\partial_{P_{t,j}} 
.\end{align}

 4. The radiation mode amplitudes  remain constant on $L^2(0,k^2)^2$, equipped with the weak topology, as $\eps \to 0$. However, this does not
describe the power transported by the radiation modes $\displaystyle \sum_{t \in \{e,o\}}\int_0^{k^2} |a_{t,\gamma}|^2 d\gamma$,  because the convergence does not hold in the strong topology of $L^2(0,k^2)^2$. 
\end{remark}

\vspace{0.05in} 
The third point in Remark \ref{rem:2} implies that the   mean mode powers satisfy the  closed system of  equations
\begin{align}
\label{eq:coupledmodeeq}
\frac{d \EE[ P_{t,j}(z)]}{dz}
=
\sum_{t' \in \{e,o\}}
\sum_{l'=1}^{N_{t'}}
\Gamma^{c}_{t,j,t',l'}   \big( \EE[P_{t',l'}(z)]-\EE[P_{t,j}(z)]\big) 
-     \Lambda^{c}_{t,j} \EE[{P_{t,j}(z)}] ,
\end{align}
with initial conditions 
\begin{equation}
P_{t,j}(0) = |a_{t,j}^{(0)}|^2, \qquad j = 1, \ldots, N_t, ~~ t \in \{e,o\}.
\end{equation}
These equations (\ref{eq:coupledmodeeq}) can be found in the literature \cite{marcuse,papakonstantinou}.
Here we  derived them from first principles, taking into account the coupling between all types of modes. 

\vspace{0.05in} \begin{remark}
\label{rem:3}
 The three terms $\{{\cal L}^j\}_{1 \le j \le 3}$ of the the limit infinitesimal generator ${\cal L}$ 
account for the mode coupling 
as follows: 

\vspace{0.05in}
1. The operator ${\cal L}^1$ accounts for  coupling of the guided modes, modeled by  the coefficients $\{\Gamma_{t,j,t',j'}^c \}_{1 \le j \le N_t, 1\leq j'\leq N_{t'}, t,t' \in \{e,o\}}.$ This coupling 
 results in the power exchange between these modes.

\vspace{0.05in} 
2. The operator ${\cal L}^2$  is due to coupling between the guided and radiation modes,  
which causes power leakage from the guided modes to the radiation ones (effective diffusion),  modeled by the coefficients 
$\{\Lambda_{t,j}^c\}_{1 \le j \le N_t, t \in \{e,o\}}$.

\vspace{0.05in}
3. The operator ${\cal L}^3$ accounts for coupling between the guided and evanescent modes, which gives the term proportional to 
$\{\kappa_{t,j}^{\rm ev}\}_{1 \le j \le N_t, t \in \{e,o\}}.$ The other term, proportional to $\{\kappa_{t,j}\}_{1 \le j \le N_t, t \in \{e,o\}}$, 
is due to the neglected terms discussed at the end of section \ref{sect:modecRand}. Both contributions  give only additional phase terms on the  
guide modes i.e.,  an effective dispersion.
\end{remark}

\section{Two weakly coupled single mode random waveguides}
\label{sect:randomTwo}
Using the results of the previous section and the decay in range  (like 
$({z}/{\eps^2})^{-2} = \eps^4/z^2$) of the contribution of the  radiation and evanescent modes,   we obtain that 
\begin{align}
{p} \Big(\frac{z}{\eps^2},x\Big) =  \sum_{t \in \{e,o\}} \sum_{j=1}^{N_t}  \frac{{a}_{t,j}(z) }{\sqrt{\beta_{t,j}}} \exp\Big( i \beta_{t,j} \frac{z}{\eps^2}\Big) 
 \phi_{t,j}(x)  +o(1) , \qquad z > 0,
 \label{eq:wavefield}
\end{align} 
where the amplitudes $\{a_{t,j}(z)\}_{1 \le j \le N_t, t \in \{e,o\}}$ are described in Theorem \ref{prop:1} and the residual $o(1)$ tends to zero as $\eps \to 0$. 

We now specialize this result 
to the case $N_e = N_o = 1$, 
where we can apply Lemma \ref{lem.singlebeta}, under the assumptions \eqref{eq:assumeD}
and \eqref{eq:dlarge},  and we can use the expansion \eqref{eq:phi_eo} of the eigenfunctions.
The goal is to derive the results stated in section \ref{sect:resRand},  in the three coupling regimes defined by the scaling relations (\ref{eq:modReg}--\ref{eq:huged}). 

\subsection{Moderate coupling regime}
\label{sect:ModTwo}
In the regime defined by the scaling relation \eqref{eq:modReg}, the expression \eqref{eq:wavefield} becomes 
\begin{align}
{p} \Big(\frac{z}{\eps^2},x\Big) =  \frac{\phi(|x|)}{\sqrt{\beta}} e^{i \beta \frac{z}{\eps^2}} \Big[ 1_{(0,\infty)}(x) u_+(z) + 
1_{(-\infty,0)}(x) u_-(z) \Big] + o(1),
 \label{eq:wavefieldTwo}
\end{align} 
with $\phi$ defined in \eqref{eq:phi} and  $\beta$ defined in Lemma \ref{lem.singlebeta}. It corresponds to 
a  single guided wave in each waveguide, with transverse profile defined by $\phi(|x|)$ and range dependent random amplitude
\begin{equation}
u_{\pm}(z) = a_{e}(z) \exp\Big( i \Delta \beta_e \frac{z}{\eps^2}\Big)  \pm a_{o}(z)  \exp\Big( i \Delta \beta_o \frac{z}{\eps^2}\Big). 
\label{eq:MC1}
\end{equation}
This is the expression of the wave field  given in section \ref{sect:resRand}, where we simplified the notation 
as $a_{t,1}(z) \leadsto a_t(z)$ and introduced
\begin{equation}
\Delta \beta_t = \beta_{t,1} - \beta = O\big( e^{-\eta d} \big), \quad t \in \{e,o\},
\label{eq:MC2}
\end{equation}
satisfying $
1 \gg |\Delta \beta_t| \gg \eps^2,
$
by equations \eqref{eq:modReg} and \eqref{eq:approxbetae}.

The statistics of the amplitudes \eqref{eq:MC1} follows from Theorem \ref{prop:1}, 
in the case $N_e = N_o = 1$. The coefficients  in the infinitesimal generator ${\cal L}$ are described in the next corollary, proved in Appendix \ref{ap:A}.

\begin{corollary}
\label{prop:2a}
The effective coefficients in the operator ${\cal L}^1$ are given by 
\begin{align}
\Gamma_{e,1,o,1}^c=\Gamma_{o,1,e,1}^c &= - \Gamma_{e,1,e,1}^c= - \Gamma_{o,1,o,1}^c = \Gamma,
\label{eq:defGammas}
\end{align}
and 
\begin{align}
\Gamma_{t,1,t',1}^1 = \Gamma,  \quad
\Gamma_{t,1,t,1}^s = 0, \qquad t,t'\in \{e,o\},
\end{align}
where $\Gamma$ has the expression \eqref{eq:defGamma}. The effective coefficients in the operator ${\cal L}^2$ are 
\begin{align}
\Lambda_{t,1}^c = \Lambda, \quad \Lambda_{t,1}^s = \Theta, \quad t \in \{e,o\},
\label{eq:defLambdaTheta}
\end{align}
where $\Lambda$ is given by \eqref{eq:defLambda} and 
\begin{align}
\Theta = 
\frac{\Delta k^4D^2}{\beta\big(\frac{2}{\eta} + D\big)} {\rm cos}^2\Big(\frac{\xi D}{2}\Big)
\int_0^{k^2} \hspace{-0.05in}\frac{d \gamma}{\pi \eta_\gamma \sqrt{\gamma}}
\Bigg[\frac{ 2\frac{\xi_\gamma^2}{\eta_\gamma^2}+ \sin^2(\xi_\gamma D) \big(1-\frac{\xi_\gamma^2}{\eta_\gamma^2}\big)}
{ 4\frac{\xi_\gamma^2}{\eta_\gamma^2}+\sin^2(\xi_\gamma D) \big(1-\frac{\xi_\gamma^2}{\eta_\gamma^2}\big)^2}\Bigg]
\nonumber \\
\times \int_0^\infty dz \, {\cal R}(z) \sin \big[(\sqrt{\gamma}-\beta)z\big]. \label{eq:defTheta}
\end{align}
The effective coeffients in the operator ${\cal L}^3$ are 
\begin{align}
\kappa_{t,1}^{\rm ev} = \kappa^{\rm ev}, \quad \kappa_{t,1} = \kappa, \quad t \in \{e,o\},
\label{eq:defkappas}
\end{align}
where 
\begin{align}
\kappa^{\rm ev} = 
\frac{\Delta k^4D^2}{\beta\big(\frac{2}{\eta} + D\big)} {\rm cos}^2\Big(\frac{\xi D}{2}\Big)
\int_{-\infty}^0 \hspace{-0.05in}\frac{d \gamma}{\pi \eta_\gamma \sqrt{|\gamma|}}
\Bigg[\frac{ 2\frac{\xi_\gamma^2}{\eta_\gamma^2}+ \sin^2(\xi_\gamma D) \big(1-\frac{\xi_\gamma^2}{\eta_\gamma^2}\big)}
{ 4\frac{\xi_\gamma^2}{\eta_\gamma^2}+\sin^2(\xi_\gamma D) \big(1-\frac{\xi_\gamma^2}{\eta_\gamma^2}\big)^2}\Bigg]
\nonumber \\
\times \int_0^\infty dz \, {\cal R}(z) \cos \big( \beta z\big) \exp(-\sqrt{|\gamma|}z). 
\end{align}
and
\begin{align}
\kappa = &  -
\frac{\Delta k^2  D^2 \xi }{\beta  \big(\frac{2}{\eta} + D\big)} {\cal R}(0)  \sin (\xi D) .
\end{align}
\end{corollary}

Using Theorem \ref{prop:1} and Corollary \ref{prop:2a}, we obtain the joint moments of the mode amplitudes 
$a_{e}(z)$ and $a_{o}(z)$:
The mean  amplitudes are given by 
\begin{align}
\EE[a_{t}(z) ] &= a_{t}^{(0)} e^{-(\Gamma+\Lambda/2) z -i (\kappa^{\rm ev} + \kappa +\Theta/2) z }, \qquad 
\end{align}
with $a_{t}^{(0)}$ defined in \eqref{eq:Amplo} for $t \in \{e,o\}$.
The mean powers and field covariance are
\begin{align}
\EE[|a_{e}(z)|^2 ] =& \frac{(|a_e^{(0)}|^2+|a_o^{(0)}|^2)}{2}
e^{-\Lambda z}+
\frac{(|a_e^{(0)} |^2-|a_o^{(0)} |^2)}{2} e^{  -(2\Gamma+\Lambda) z },\label{eq:aesquared}\\
\EE[|a_{o}(z)|^2 ] =& \frac{(|a_e^{(0)}|^2+|a_o^{(0)}|^2)}{2}
e^{-\Lambda z}-
\frac{(|a_e^{(0)} |^2-|a_o^{(0)} |^2)}{2} e^{  -(2\Gamma+\Lambda) z },\label{eq:aosquared}\\
 \EE [ a_{o}(z)\overline{a_{e}(z)}] =&
 (a_o^{(0)}\overline{a_e^{(0)}})e^{  -(\Gamma+\Lambda) z}, \label{eq:aeaocovar}
\end{align}
and the  fourth-order moments are
\begin{align}
\nonumber
& \begin{pmatrix}
\EE[|a_{e}(z)|^4 ] \\
\EE[|a_{o}(z)|^4 ] \\
2\EE[|a_{e}(z)|^2|a_{o}(z)|^2 ]
\end{pmatrix}
=
\frac{(|a_e^{(0)}|^2+|a_o^{(0)}|^2)^2)}{3} e^{-2\Lambda z} \begin{pmatrix} 1\\1\\1\end{pmatrix} \\
\nonumber
& \hspace{1in}-
\frac{(|a_e^{(0)}|^4-|a_o^{(0)}|^4)}{2} e^{ -(2\Lambda +2\Gamma) z } \begin{pmatrix} -1\\1\\0\end{pmatrix}
\\
&\hspace{1in}  -
\frac{(|a_e^{(0)}|^4+|a_o^{(0)}|^4-4 |a_e^{(0)}|^2|a_o^{(0)}|^2)}{6} e^{ -(2\Lambda +6\Gamma) z} \begin{pmatrix} -1\\-1\\2\end{pmatrix}
. \label{eq:fourthmoment}
\end{align} 

We conclude from  equations (\ref{eq:aesquared}--\ref{eq:aosquared})  and \eqref{eq:fourthmoment} that 
\begin{align}
\EE\big[ |a_{e}(z)|^2+|a_{o}(z)|^2 \big] =& ( |a_e^{(0)}|^2+|a_o^{(0)}|^2) e^{-\Lambda z},\\
{\rm Var}\big( |a_{e}(z)|^2+|a_{o}(z)|^2\big) =& 0.
\end{align}
This implies the results \eqref{eq:detpower} and \eqref{eq:powerdecay}, which model the effective leakage of power from the guided modes 
to  the radiation modes. 
Equations (\ref{eq:aesquared}--\ref{eq:aeaocovar}) also give the expression \eqref{eq:powertransf} of the imbalance of power 
between the waveguides. 

\begin{remark}
\label{rem:10}
These results show  that of all the coefficients in Corollary \ref{prop:2a},  $\Gamma$ and $\Lambda$ are the important ones.
The coefficient $\Gamma$ determines the coupling between the two guided modes,
whereas $\Lambda$  determines the power leakage to the radiation modes. Note from \eqref{eq:defLambda}--\eqref{eq:defGamma} that 
$\Gamma$ depends on power spectral density $\widehat {\cal R}$, the Fourier transform of the covariance ${\cal R}$, evaluated  at  zero wavenumber,
while $\Lambda$ depends on  $\widehat {\cal R}$ evaluated at  wavenumbers larger than $\beta-k$, where we recall that 
$\beta \in (k,n k)$.  This result is consistent with the observation in  \cite{marcuse69} that radiation losses should only depend on the power spectral density evaluated at  wavenumbers between $\beta-k$ and $\beta+k$. 
\end{remark}

It is possible to encounter situations where $\Gamma > 0$ and $\Lambda = 0$, or the opposite.
Indeed, if the power spectral density  vanishes for wavenumbers larger than $\beta-k$, then 
$\Lambda = 0$. Otherwise, if $\widehat {\cal R}(0) = 0$  (this happens for instance
when the random processes $\nu_j$ are derivatives of stationary processes), then $\Gamma = 0$.
However, in general, both $\Gamma$ and $\Lambda$ are positive.
\subsection{Weak coupling regime}
\label{sect:twoWeak}
In the regime defined by the scaling relation \eqref{eq:larged}, the wave field has the same expression 
\eqref{eq:wavefieldTwo}, but the amplitudes $u_\pm(z)$ have different statistics. They are obtained as the 
$\eps \to 0$ limit of 
\begin{equation}
u_\pm^{(\eps)}(z) = a_{e,1}^{(\eps)}(z) e^{i \beta' \theta z} \pm a_{o,1}^{(\eps)}(z) e^{-i \beta' \theta z},
\label{eq:Wk1}
\end{equation}
where we used the expression of $\beta_{e,1}$ and $\beta_{o,1}$ given in equations  (\ref{eq:approxbetae}) and
set $\exp(-\eta d) =\eps^2 \theta$ with $\theta \in [0,\infty)$.
To obtain this limit, we let 
\begin{equation}
\alpha_e^{(\eps)}(z) = a_{e,1}^{(\eps)}(z) e^{i \beta' \theta z}, \qquad \alpha_o^{(\eps)}(z) = a_{o,1}^{(\eps)}(z) e^{-i \beta' \theta z},
\label{eq:Wk2}
\end{equation}
and then use equations \eqref{eq:ForwEq1}--\eqref{eq:ForwEq2} to derive a closed system of stochastic 
differential equations satisfied by the process $\big( \alpha_t^{(\eps)}(z), \big(a_{t,\gamma}^{(\eps)}(z)\big)_{\gamma \in (0,k^2)}, t \in \{e,o\}\big)$.
The limit  of this process is given in the next theorem, with proof outlined  in Appendix \ref{ap:Theom2}. As was the case in 
Theorem \ref{prop:1}, the infinitesimal generator does not involve the radiation mode amplitudes, 
so we  describe directly the limit of $\big( (\alpha_t^{(\eps)}(z), t \in \{e,o\}\big)$.

\vspace{0.05in}
\begin{theorem}
\label{thm:2}
In the limit $\eps \to 0$, the process  $\big({\alpha}_{e}^{(\eps)} (z),{\alpha}_{o}^{(\eps)}(z) \big)$
converges in distribution in ${\cal C}^0\big([0,{L}], \CC^{2} \big)$ to
the Markov process  $\balpha(z) = (\alpha_{e} (z),\alpha_{o}(z))$
with infinitesimal generator
\begin{align}
 {\cal L}_{\balpha,\theta} = {\cal L}_\balpha^1+ {\cal L}_\balpha^2+ {\cal L}_\balpha^3+
\theta {\cal L}_\balpha^4, \label{eq:tildeLtheta}
\end{align}
where $\theta$ is defined in \eqref{eq:weak3},  and the operators are defined by 
\begin{align}
\nonumber
 {\cal L}_\balpha^1 = & 
 \frac{\Gamma}{2}  \Big( 
 2 (|\alpha_o|^2+|\alpha_e|^2) (\partial_{{\alpha_o}} \partial_{\overline{\alpha_{o}}}+\partial_{{\alpha_e}} \partial_{\overline{\alpha_{e}}})
+2 (\alpha_o\overline{\alpha_e}+\overline{\alpha_o}\alpha_e) (\partial_{{\alpha_o}} \partial_{\overline{\alpha_{e}}}+\partial_{{\alpha_e}} \partial_{\overline{\alpha_{o}}})
  \\
\nonumber
&
 -
(\alpha_{o}^2+\alpha_e^2) ( \partial_{\alpha_o}^2+ \partial_{\alpha_e}^2)
-
(\overline{\alpha_o}^2+  \overline{\alpha_e}^2)(\partial_{\overline{\alpha_o}}^2+\partial_{\overline{\alpha_{e}}}^2)
 -  4 \alpha_o \alpha_e  \partial_{\alpha_o} \partial_{\alpha_e}
\\
 &  - 4 \overline{\alpha_o} \overline{\alpha_e} \partial_{\overline{\alpha_o}} \partial_{\overline{\alpha_{e}}}
 - 
2(\alpha_o\partial_{\alpha_o}+\alpha_{e}  \partial_{\alpha_e})
-
2(\overline{\alpha_o} \partial_{\overline{\alpha_o}}+  \overline{\alpha_e}\partial_{\overline{\alpha_{e}}})
 \Big),
\label{def:tildeL1}
\end{align}
and 
\begin{align}
{\cal L}_\balpha^2
=&
-\frac{1}{2} \sum_{t\in \{e,o\}}
( \Lambda  +i \Theta ) {\alpha}_{t} \partial_{{\alpha}_{t}} 
+
( \Lambda - i\Theta )  \overline{{\alpha}_{t}} \partial_{\overline{{\alpha}_{t}}}  ,
 \label{eq:tildeL2}
 \\
 {\cal L}_\balpha^3
=&
i  \big( \kappa + \kappa^{\rm ev} \big) \sum_{t\in \{e,o\}}
   \big( {\alpha}_{t} \partial_{{\alpha}_{t}} 
- \overline{{\alpha}_{t}} \partial_{\overline{{\alpha}_{t}}} \big) , \label{eq:tildeL3}
\\
{\cal L}_\balpha^4
=& i 
  \beta'    \big(  \alpha_{e} \partial_{\alpha_{e}} - \overline{\alpha_{e}} \partial_{\overline{\alpha_{e}}}
-
\alpha_{o} \partial_{\alpha_{o}} + \overline{\alpha_{o}} \partial_{\overline{\alpha_{o}}} \big)  ,
\label{def:tildecalLs4}
\end{align}
with coefficients given in Corollary \ref{prop:2a}.
\end{theorem}

The result of Theorem \ref{thm:2} can be reformulated as stated in the following corollary, by the change of variables:
\[u_\pm^{(\eps)}(z) = \alpha_{e}^{(\eps)}(z)\pm \alpha_{o}^{(\eps)}(z).\]

\begin{corollary}
In the limit $\eps \to 0$, the process
$(u_+^{(\eps)}(z) ,u_-^{(\eps)}(z) )$ converges in distribution in ${\cal C}^0([0,L], \CC^2)$ to the Markov process ${\itbf u}(z) = (u_+ (z), u_-(z) )$ with infinitesimal generator
$
{\cal L}_{{\itbf u},\theta} = {\cal L}_{u_+}+{\cal L}_{u_-}+ \theta {\cal L}_{u_+,u_-},
$
defined by 
\begin{align}
\nonumber
{\cal L}_{u} = & \Gamma \big( 2 |{u}|^2 \partial_{u} \partial_{\overline{{u}}}
-   {u}^2 \partial_{u}^2
-\overline{{u}}^2 \partial_{\overline{{u}}}^2
-    {u}\partial_{u} 
-  \overline{{u}} \partial_{\overline{{u}}}\big)\\
&- \frac{ \Lambda + i\Theta }{2} {u}\partial_{u} 
 - \frac{  \Lambda - i\Theta }{2} \overline{{u}} \partial_{\overline{{u}}} 
 + i (\kappa +\kappa^{\rm ev})  \big({u} \partial_{u}- 
\overline{{u}} \partial_{\overline{{u}}}\big) ,
\label{def:calLu}
\\
{\cal L}_{u_+,u_-} =& i  \beta' \big( u_+ \partial_{u_-} + u_- \partial_{u_+}
- \overline{u_+}\partial_{\overline{u_-}} - \overline{u_-}\partial_{\overline{u_+}} \big).
 \end{align}
 \end{corollary}

Note that for any test function $F$, the generator gives
 \begin{equation}
{\cal L}_{{\itbf u},\theta}  F(|u_+|^2+|u_-|^2) = -  \Lambda (|u_+|^2+|u_-|^2) F'(|u_+|^2+|u_-|^2).
 \label{eq:powerDecay1}
 \end{equation}
This shows that $|u_+(z)|^2+|u_-(z)|^2 = 2 (|\alpha_o(z)|^2+|\alpha_e(z)|^2)$ is deterministic 
and satisfies \eqref{eq:powerdecay}, as stated in 
 section \ref{sect:resWeak}. The total power transported by the two guided modes decays as $\exp (- \Lambda z) $
with probability one, due to an effective leakage towards the radiation modes:
\begin{equation}
|u_+(z)|^2+|u_-(z)|^2 = 2 (|a_o^{(0)}|^2+|a_e^{(0)}|^2) \exp(-\Lambda z) .
\end{equation}

To determine the  imbalance of power between the two guided modes
\begin{align}
 \mathcal{P} (z)   =\frac{|u_+(z)|^2 - |u_-(z)|^2}{ |u_+(z)|^2+|u_-(z)|^2} ,
\end{align}
we apply the infinitesimal generator ${\cal L}_{{\itbf u},\theta} $ to the two functions 
\[(u_+,\overline{u_+},u_-,\overline{u_-})\mapsto |u_+|^2-|u_-|^2, \qquad 
(u_+,\overline{u_+},u_-,\overline{u_-}) \mapsto \overline{u_+} u_-,
\]
and we get
\begin{align} 
{\cal L}_{{\itbf u},\theta}  (|u_+|^2-|u_-|^2) &= -  \Lambda (|u_+|^2-|u_-|^2)  
 -4\theta \beta' {\rm Im}(\overline{u_+} u_-)   ,\\
{\cal L}_{{\itbf u},\theta}  (\overline{u_+} u_-) &= (-2\Gamma-  \Lambda)(\overline{u_+} u_-)  
 +i \theta \beta' (|u_+|^2-|u_-|^2). \label{eq:ImbalDeriv}
 \end{align}
 Denoting ${\cal I}(z) = 2 {\rm Im}(\overline{u_+} u_-(z)) / (|u_+(z)|^2+|u_-(z)|^2)$, 
we find that 
\begin{align*}
&
\partial_z \EE[ {\cal P}(z) ]= -2\theta \beta' \EE[ {\cal I}(z)  ],\\
&
\partial_z \EE[ {\cal I}(z) ]= -2\Gamma \EE[ {\cal I}(z) ]  + 2\theta\beta'\EE[ {\cal P}(z) ],
\end{align*}
which gives the harmonic oscillator equation \eqref{eq:harmosc} satisfied by $\EE[  \mathcal{P} (z) ] $.

\subsection{Very weak coupling regime}
When the separation  $d$ between the waveguides is so large  that 
$\exp( -\eta d) \ll \eps^2$, we obtain from \eqref{eq:approxbetae}
that 
\begin{equation}
\beta_{t,1} = \beta + o(\eps^2), \qquad t \in \{e,o\}.
\end{equation}
The $\eps \to 0$ limit  in this case is similar to that in the previous section, and the result is the same as setting 
$\beta' = 0$ in  Theorem \ref{thm:2}. The expression \eqref{eq:wavefieldTwo} of the wave field still holds, 
but the mode amplitudes  have different statistics, as described in the next corollary:

\vspace{0.05in}
\begin{corollary}
\label{prop:2c}
In the limit $\eps \to 0$, the process
 $( u_+^{(\eps)} (z) , u_-^{(\eps)}(z) )$ 
converges in distribution in ${\cal C}^0([0,{L}], \CC^{2} )$ to the process $(u_+(z),u_-(z))$ where $u_+(z)$ and $u_-(z)$ are independent and identically distributed Markov processes with infinitesimal generator ${\cal L}_u$ defined by (\ref{def:calLu}).
\end{corollary}

 Note that, for any test function $F$, we have
 $$
 {\cal L}_u F(|u|^2) = - \Lambda  |u|^2 F'(|u|^2) .
 $$
This shows that the processes $(|u_+(z)|^2,|u_-(z)|^2)$ are deterministic and exponentially decaying with the rate $\Lambda$.
We conclude  that 
$(|u_+^{(\eps)}|^2(z), |u_-^{(\eps)}|^2(z))$ 
converges in probability to the deterministic function  $(|u_+(z)|^2,|u_-(z)|^2)$ 
that decays exponentially as $\exp ( - \Lambda z )$. Therefore, there is no mode coupling in this regime,  except 
between guided and radiation modes, which results in effective leakage. 
In particular, the imbalance of power between the two waveguides is constant, as stated in section 
\ref{sect:resVWeak}.

\section{Summary}
\label{sect:sum}
In this paper we introduced an analysis of wave propagation in a directional coupler 
consisting of two parallel step-index waveguides. The waveguide effect is due to a 
medium of high index of refraction separated from a uniform background by randomly 
fluctuating interfaces. The fluctuations occur on a length scale (correlation length) that is similar 
to the wavelength and have small amplitude modeled by a dimensionless parameter $\eps$ satisfying 
$0 < \eps \ll 1$.  The analysis is based on the decomposition of the wave field in a complete set 
of guided, radiation and evanescent modes. It accounts for the interaction of all these modes and derives a closed system of 
stochastic differential equations for the guided and radiation mode amplitudes. These equations are driven by 
the random fluctuations of the interfaces and model the net scattering effect that becomes significant at distances 
of propagation of the order of $\eps^{-2}$.  We analyze them in the asymptotic limit $\eps \to 0$, under the assumption that the 
covariance function of the fluctuations is smooth. This allows us to use the forward scattering approximation and 
characterize the limit mode amplitudes as a Markovian process with infinitesimal generator that we calculate explicitly. 
The analysis applies to waveguides that support an arbitrary number of guided modes. Since many 
directional couplers use  single guided mode waveguides, we studied in detail this case and obtained 
a detailed  quantification of the transfer of power for three different regimes, where the coupling 
between the guided waves is stronger or weaker, depending on how far apart the waveguides are.  In all regimes 
there is a self-averaging power leakage from the guided to the radiation modes.  This is the only coupling if the waveguides are very far apart. Otherwise, the guided modes are coupled and the random boundary fluctuations induce a blurring of the periodic transfer of power that would occur
in the absence of the random fluctuations. We quantified this blurring and showed that at sufficiently long distances 
of propagation the power becomes evenly distributed among the waveguides, independent of the initial condition defined by the wave source.

\section*{Acknowledgements}
This research is supported in part by the Air Force Office of Scientific Research under award number FA9550-18-1-0131 and in part
by the NSF grant DMS1510429.
\appendix
\section{Proof of Lemma \ref{lem.singlebeta}}
\label{sect:ap0}
The solution $\beta$ of \eqref{eq:barbeta1} can be obtained as follows: First, note that $cos(q)/q$ is a monotone,
strictly decreasing function in the interval $q \in (0,\pi/2)$, with
\[\lim_{q\searrow  0}  \frac{\cos q}{q} = \infty, \qquad \lim_{q \nearrow \pi/2} \frac{\cos q}{q} = 0.\]
Thus, there is a unique solution $q^\star$ of
\[
\frac{\cos q}{q} = \frac{2}{k D \sqrt{n^2-1}},
\]
which gives
\[
\tan q^\star = \frac{\Big[1-\Big(\frac{2 q^\star}{ k D \sqrt{n^2-1}}\Big)^2\Big]^{1/2}}{\frac{2q^\star}{k D \sqrt{n^2-1}}}.
\]
Substituting
$
q^\star = \frac{D}{2} \sqrt{n^2 k^2 - \beta^2}
$
into  this equation we obtain that \eqref{eq:barbeta1} is satisfied, with $
\beta = \sqrt{k^2 n^2 - \Big(\frac{2 q^\star}{D}\Big)^2}.
$
Obviously,  $\beta < k n$ and the condition $\beta > k$ is consistent with
\[
q^\star < \frac{k D}{2}\sqrt{n^2-1} < \frac{\pi}{2},
\]
where the last inequality is by the assumption \eqref{eq:assumeD}.

A similar construction shows that \eqref{eq:barbeta2} has no solution in the interval $ (k,nk)$.
Therefore, the right-hand sides of equations \eqref{eq:dispersione} and \eqref{eq:dispersiono}
vanish at the single root $\beta \in (k,nk)$ of \eqref{eq:barbeta1}.

If  $d$ is large, then the left-hand sides of equations \eqref{eq:dispersione} and \eqref{eq:dispersiono}
are small, and  by the implicit function theorem, they have unique solutions, close to  $\beta$ calculated above.
Consequently, $N_e=N_o=1$.

\section{Proof of Corollary \ref{prop:2a}}
\label{ap:A}

Theorem \ref{prop:1} states that 
\begin{equation}
\Gamma^{c}_{e,1,o,1}=
\frac{\Delta k^4}{2 \beta_{e}\beta_{o}}
\int_0^\infty \EE\big[ C_{e,1,o,1}(0) C_{e,1,o,1}(z) \big] \cos \big[ (\beta_{e}-\beta_{o})z\big] dz,
\label{eq:Gammaeo}
\end{equation}
where $\beta_e \approx \beta_o \approx \beta$ and the expectation follows by definitions \eqref{expres:Cttp}, \eqref{eq:phi}
and equation \eqref{eq:phi_eo},
\begin{align*}
\EE\big[ C_{e,1,o,1}(0) C_{e,1,o,1}(z) \big] 
&\approx 2 {\cal R}(z) D^2\Big[ \phi_{e,1}^2\Big(\frac{d}{2}\Big) \phi_{o,1}^2 \Big(\frac{d}{2}\Big) + \phi_{e,1}^2\Big(\frac{d}{2}+D\Big)
 \phi_{o,1}^2 \Big(\frac{d}{2}+D\Big)\Big], \\
 &\approx \frac{4 D^2}{\Big(\frac{2}{\eta} + D \Big)^2} {\cal R}(z) \cos^4 \Big(\xi\frac{D}{2}\Big).
\end{align*}
Here the approximation is up to an error of the order of $\exp(-\eta d) \ll 1$, which we neglect in this appendix.
Substituting in \eqref{eq:Gammaeo} we get the results \eqref{eq:defGammas} and \eqref{eq:defGamma}.
Proceeding similarly, we find $\Gamma^1_{t,1,t',1}=\Gamma,$ for all $ t,t'\in \{e,o\}.$

The expression of $\Lambda^{c}_{t,1}$ in Theorem \ref{prop:1} is 
\begin{equation}
\Lambda_{t,1}^{c} =  \sum_{t'\in \{e,o\}}  \int_0^{k^2} d \gamma  \, \frac{\Delta k^4}{2\sqrt{\gamma}\beta_{t,1}}
\int_0^\infty dz \, \EE \big[ C_{t,1,t',\gamma}(0) C_{t,1,t',\gamma}(z) \big] 
\cos \big( (\sqrt{\gamma}-\beta_{t,1})z\big)  ,
\label{eq:calcLambda}
\end{equation}
where the expectation is calculated using the definition \eqref{eq:defCoef2},
\begin{equation}
\label{eq:expectationCga}
\EE \big[ C_{t,1,t',\gamma}(0) C_{t,1,t',\gamma}(z) \big] 
=2 {\cal R}(z) D^2\Big[ \phi_{t,1}^2\Big(\frac{d}{2}\Big) \phi_{t',\gamma}^2 \Big(\frac{d}{2}\Big) + 
\phi_{t,1}^2\Big(\frac{d}{2}+D\Big) \phi_{t',\gamma}^2 \Big(\frac{d}{2}+D\Big)\Big] .
\end{equation}
As explained above, $\phi_t^2 \approx \phi^2$ for $t \in \{e,o\}$,  and definition \eqref{eq:phi} gives that 
\[
\phi^2\Big(\frac{d}{2}\Big) = \phi^2 \Big(\frac{d}{2}+D\Big) = \Big(\frac{2}{\eta} + D \Big)^{-1} \cos^2 \Big(\xi \frac{D}{2} \Big).
\]
We also have from  the expressions of $\phi_{e,\gamma}$ and $\phi_{o,\gamma}$ in section \ref{sect:continSpect} that 
\begin{align*}
 \phi_{e,\gamma}^2 \Big(\frac{d}{2}\Big) + \phi_{e,\gamma}^2 \Big(\frac{d}{2}+D\Big) 
&= A_{e,\gamma}^2 \Big\{\frac{\xi_{\gamma}^2}{\eta_{\gamma}^2} \cos^2 (\eta_{\gamma}\frac{d}{2})\\
&+\big[ -\sin (\xi_{\gamma} D) \sin (\eta_{\gamma} \frac{d}{2}) +\frac{\xi_{\gamma}}{\eta_{\gamma}}
\cos(\xi_{\gamma} D) \cos(\eta_{\gamma} \frac{d}{2}) \big]^2 \Big\} ,
\end{align*}
with $A_{e,\gamma}$ given by (\ref{def:Aegamma}).  Substituting in \eqref{eq:expectationCga}, we obtain 
\begin{align*}
& \EE \big[ C_{t,1,e,\gamma}(0) C_{t,1,e,\gamma}(z) \big] 
=
\frac{D^2}{\pi \eta_\gamma \Big(\frac{2}{\eta} + D\Big)} {\cal R}(z) \cos^2 \big(\frac{\xi D}{2} \big)
  \\
&\times \frac
{
\cos (\eta_{\gamma} d) \Big[2\frac{\xi_{\gamma}^2}{\eta_{\gamma}^2}-\sin^2(\xi_{\gamma}D) (1+\frac{\xi_{\gamma}^2}{\eta_{\gamma}^2})\Big]
-\sin(\eta_{\gamma} d) \frac{\xi_{\gamma}}{\eta_{\gamma}} \sin(2 \xi_{\gamma} D) + 
2\frac{\xi_{\gamma}^2}{\eta_{\gamma}^2}+\sin^2(\xi_{\gamma}D) (1-\frac{\xi_{\gamma}^2}{\eta_{\gamma}^2})
}
{
\cos (\eta_{\gamma} d)  \sin^2(\xi_{\gamma} D) \Big(\frac{\xi_{\gamma}^4}{\eta_{\gamma}^4}-1\Big)
+ \sin(\eta_{\gamma} d) \frac{\xi_{\gamma}}{\eta_{\gamma}} \Big(\frac{\xi_{\gamma}^2}{\eta_{\gamma}^2}-1\Big)  \sin (2 \xi_{\gamma} D) 
+ 2\frac{\xi_{\gamma}^2}{\eta_{\gamma}^2}+ \sin^2(\xi_{\gamma}D) (1-\frac{\xi_{\gamma}^2}{\eta_{\gamma}^2})^2
} .
\end{align*}
As $d$ is large, we take the weak limit $d \to+\infty$ in this expression, seen as a function of $\gamma$, and we get 
$$
 \EE \big[ C_{t,1,e,\gamma}(0) C_{t,1,e,\gamma}(z) \big] 
= \frac{D^2}{\pi \eta_\gamma \Big(\frac{2}{\eta} + D\Big)} {\cal R}(z)  \cos^2 \Big(\frac{\xi D}{2} \Big)
\Bigg[ \frac{ 4\frac{\xi_{\gamma}^2}{\eta_{\gamma}^2}+2\sin^2(\xi_{\gamma} D) \big(1-\frac{\xi_{\gamma}^2}{\eta_{\gamma}^2}\big)}
{ 4\frac{\xi_{\gamma}^2}{\eta_{\gamma}^2}+\sin^2(\xi_{\gamma} D) \big(1-\frac{\xi_{\gamma}^2}{\eta_{\gamma}^2}\big)^2}\Bigg]
$$
 by the following calculation:
for any $\alpha_1,\alpha_2,\alpha_3,\alpha_4,\alpha_5,\alpha_6$, with $\alpha_4 > \sqrt{\alpha_5^2+\alpha_6^2}>0$, we have
$$
\frac{1}{2\pi} \int_0^{2\pi} \frac{\alpha_1 +\alpha_2\cos(s) +\alpha_3 \sin(s) }{\alpha_4+\alpha_5 \cos (s) +\alpha_6 \sin(s)} ds
=
 \frac{\alpha_1 +\frac{\alpha_2\alpha_5+\alpha_3\alpha_6}{\alpha_5^2+\alpha_6^2} \big( \sqrt{\alpha_4^2-\alpha_5^2-\alpha_6^2} -\alpha_4 \big) }{\sqrt{\alpha_4^2-\alpha_5^2-\alpha_6^2}}
.
$$
This gives
$\Lambda_{e,1}^{c} = \Lambda$, as stated in Corollary \ref{prop:2a}.
We can deal with 
\begin{align*}
 \phi_{o,\gamma}^2 \Big(\frac{d}{2}\Big) + \phi_{o,\gamma}^2 \Big(\frac{d}{2}+D\Big) 
=& A_{o,\gamma}^2 \Big\{ \frac{\xi_{\gamma}^2}{\eta_{\gamma}^2} \sin^2 \Big(\eta_{\gamma}\frac{d}{2}\Big)\\
&+\Big[ \sin (\xi_{\gamma} D) \cos (\eta_{\gamma} \frac{d}{2}) +\frac{\xi_{\gamma}}{\eta_{\gamma}}
\cos(\xi_{\gamma} D) \sin\Big(\eta_{\gamma} \frac{d}{2}\Big) \Big]^2 \Big\}
\end{align*}
and $\EE \big[ C_{t,1,o,\gamma}(0) C_{t,1,o,\gamma}(z) \big] $
in the same way, which gives
$\Lambda_{o,1}^{c} = \Lambda$.
The coefficients $\Lambda_{t,1}^{s}$, $\kappa_{t,1}$, $t\in \{e,o\}$ are obtained in a similar way.

\section{Proof of Theorem \ref{thm:2}}
\label{ap:Theom2}
The result follows from an extended version of Theorem \ref{prop:1}, when all $\beta_{t,j}$ are equal, 
as was considered in \cite{kohler77}.
We obtain
\begin{align*}
{\cal L}_{\balpha,\theta} = & \hspace{-0.2in}
\sum_{t_1,t_2,t_3,t_4\in \{e,o\}} \frac{\Delta k^4}{4\beta^2} \int_0^\infty \EE \big[ C_{t_1,1,t_2,1}(z)C_{t_3,1,t_4,1}(0) \big] dz 
\big(
\alpha_{t_2}\partial_{\alpha_{t_1}} \overline{\alpha_{t_4}}\partial_{\overline{\alpha_{t_3}}} 
-
\alpha_{t_2}\partial_{\alpha_{t_1}} \alpha_{t_4}\partial_{\alpha_{t_3}} 
\big)\\
&-\hspace{-0.2in}
\sum_{t_1,t_2,t'\in \{e,o\}} 
\int_0^{k^2} \frac{\Delta k^4}{4\beta \sqrt{ \gamma'}} 
 \int_0^\infty \EE \big[ C_{t',\gamma',t_2,1}(z)C_{t_1 ,1,t',\gamma'}(0) \big]
e^{i (\beta -\sqrt{\gamma'})z} dz d\gamma'
\alpha_{t_2}\partial_{\alpha_{t_1}}   \\
&
+i(\kappa+\kappa^{\rm ev})  \big(  {\alpha}_{e} \partial_{{\alpha}_{e}} +
{\alpha}_{o} \partial_{{\alpha}_{o}}\big)
+ i\theta \beta' \big(  {\alpha}_{e} \partial_{{\alpha}_{e}} -
{\alpha}_{o} \partial_{{\alpha}_{o}}\big) + c.c.,
\end{align*}
where $c.c.$ stands for complexe conjugate.
The statement of Theorem \ref{thm:2} follows from this expression, by using the definition (\ref{expres:Cttp}) of the random coefficients $
C_{t_1,1,t_2,1}(z)$ and $C_{t',\gamma',t_2,1}(z)$ and  the symmetry of the eigenfunctions i.e.,  that $\phi_{e,1}$ is even and $\phi_{o,1}$ is odd.


\end{document}